\documentclass[11pt, a4paper]{article}

\usepackage[
top = 1.2in,
bottom=1.2in,
left=1.35in,
right=1.35in
]{geometry}
\usepackage{amsmath, amssymb, amsthm}
\usepackage{mathtools}
\usepackage{float}
\usepackage{graphicx}
\usepackage{mathrsfs} 
\usepackage{booktabs} 
\usepackage{subcaption} 

\usepackage{microtype}
\numberwithin{equation}{section} 

\usepackage{hyperref}
\hypersetup{colorlinks=true, citecolor = blue}

\theoremstyle{plain}
\newtheorem{theorem}{Theorem}
\newtheorem{lemma}{Lemma}
\newtheorem{proposition}[theorem]{Proposition}

\theoremstyle{definition}
\newtheorem{example}{Example}
\newtheorem{remark}{Remark}

\newcommand{\spr}[6]{\textsc{#1}: #2. \textit{#3} {\textbf{#4}}, #5 (#6) }

\newcommand{\sprthesis}[6]{\textsc{#1}: #2. \textit{#3}. #4, #5, #6.}

\newcommand{\sprbook}[5]{\textsc{#1}: #2. \textit{#3}, #4, #5.}

\newcommand{\sprconf}[7]{\textsc{#1}: #2. \textit{#3}. #4. #5, #6, #7.}

\newcommand{\sprbookchp}[7]{\textsc{#1}: #2, \textit{#3}. #4. #5, #6, #7.}

\usepackage{tikz}

\tikzset{
  style1/.style={black, line width=1.0pt, line cap=round,
    dash pattern=on 0.1pt off 5.1pt},                       
  style2/.style={black, line width=1.0pt,
    dash pattern=on 5.4pt off 3.7pt},                       
style3/.style={black, line width=1.0pt, line cap=round,
  dash pattern=on 4pt off 3pt on 0.5pt off 3pt},      
  style4/.style={black, line width=1.0pt, line cap=round,
    dash pattern=on 5.8pt off 2.7pt on 0.1pt off 3.1pt
                  on 0.1pt off 4.1pt},                      
  style5/.style={black, line width=1.0pt},                     
}

\DeclareRobustCommand{\linestyle}[2][0.8cm]{%
  \tikz[baseline=-0.5ex]{\draw[#2] (0,0) -- (#1,0);}%
}

\newcommand{\bm}[1]{\boldsymbol{#1}}
\newcommand{\Px}{\mathbf{P}}
\newcommand{\Ex}{\mathbf{E}}
\newcommand{\dif}{\mathrm{d}} 
\newcommand{\Dx}{\mbox{$\frac{\partial}{\partial x}$}}
\DeclareMathOperator{\var}{Var}
\DeclareMathOperator{\dHe}{He}

\newcommand{\domSpace}{\mathscr{X}}
\newcommand{\rangeSpace}{\mathscr{Y}}
\newcommand{\domSubSpace}{\mathscr{Z}}
\newcommand{\contSpace}{\mathscr{C}}
\newcommand{\contSpaceZero}{\mathscr{U}}
\newcommand{\domain}{\mathscr{U}}
\newcommand{\cmSpace}{\mathscr{H}}
\newcommand{\LpSpace}{\mathscr{L}}
\DeclareMathOperator{\LinSpace}{Lin}

\title{Taylor expansions for first-passage time \\
distributions of Brownian motion \\
with moving boundaries}
\date{}
\author{Vincent Liang\footnote{Department of Mathematics, City University of Hong Kong, Kowloon, Hong Kong; e-mail: vliang@cityu.edu.hk}}

\begin{document}

\maketitle

\begin{abstract}
We analyse the distribution of the first-passage time of Brownian motion to a time-dependent boundary as a function of the boundary. We prove that the boundary-to-distribution map is infinitely many times differentiable along the Cameron--Martin space~$\mathscr{H},$ and obtain a with a bound on the truncation error Taylor expansion with a bound on the truncation error. We obtain identities connecting (i) the $\mathscr{H}$-derivatives of the first-passage time distribution, (ii) the Wiener--It\^o chaos kernels of the indicator that the Brownian motion stays below the boundary, and (iii) Krylov--Veretennikov--Zvonkin-type operators. These identities yield relatively tractable representations for the $\mathscr{H}$-derivatives of the first-passage time distribution when the taboo transition density is available, notably for linear boundaries. As an application of the above results we propose Taylor approximations up to second order of the first-passage time distribution. In the examples considered, the first-order approximation is more accurate than the classical tangent approximation and its hazard-rate variant, while remaining straightforward to implement. 
\medskip 


\medskip 

\end{abstract}

\section{Introduction}

The moving boundary crossing problem for Brownian motion is of fundamental interest in probability due to the universality of the Brownian motion as a limit of sequences of scaled random walks~\cite{Kolmogorov1931_OMEN,Kolmogorov1933b}. The problem appears in a variety of contexts, including probability theory~\cite{Denisov2018}, statistical inference~\cite{Doob1949,Kolmogorov1933a}, sequential analysis~\cite{Lerche1986,Siegmund1985}, optimal stopping~\cite{Peskir2006}, neuroscience~\cite{Sacerdote2013} and mathematical finance~\cite{Kunitomo1992}.  

Despite the long history of the problem, cases of moving boundaries for which the first-passage time distribution is known in closed-form are largely confined to boundaries obtained using the method of images~\cite{Kahale2008,Lerche1986} and to a certain family of transformations of those boundaries~\cite{Alili2014}. For boundaries outside of this class, approximation methods are typically required.  Although numerical methods such (such as those based on integral equations~\cite{Peskir2002}) are effective for approximately computing first-passage time distributions, the analytical properties of these distributions are of independent interest since first-passage times of moving boundaries can arise as auxiliary objects in applications. For example, they have recently been shown to arise in the study of hydrodynamic
limits of particle systems~\cite{Klump2023} and Volterra processes~\cite{AbiJaber2026}. 

Because closed-form expressions for first-passage time distributions are not available outside the narrow class described above, a natural alternative is to preserve analytical tractability by using the well-known \emph{boundary approximation method}. In this approach, the original boundary~$g$ is approximated with another close boundary~$f$ whose first-passage time distribution is simpler to analyse or compute, and the approximation error is then quantified. Such approximations are used in both the \emph{direct} first-passage time problem~\cite{Borovkov2005, Wang1997}, which asks for the first-passage time distribution for a given boundary, and the \emph{inverse} first-passage time problem~\cite{Gur2021,Zucca2009}, which asks for the boundary which yields a given first-passage time distribution. A popular choice for the approximating boundary~$f$ is the linear boundary which is tangent to the original boundary~$g$ at the time point~$t$ of interest. The density of the first-passage time of the boundary~$f$ is given by the well-known Bachelier--L\'evy formula (see e.g.~\cite{Lerche1986}). This is known as the \emph{tangent approximation} method, which was introduced in~\cite{Strassen1967} and independently in~\cite{Daniels1974}. Another common choice for~$f$ is a piecewise-linear approximation of~$g$. This approximation is useful because its non-crossing probability can be expressed as a multiple Gaussian integral (see e.g.\ \cite{Novikov1999,Potzelberger2001,Wang1997}), and because the approximation can be refined to approximate a continuous boundary~$g$ to arbitrary accuracy (see e.g.\ its use in~\cite{Funaki2007,Liang2023}).

A necessary step in the boundary approximation method is to bound  the error that arises in the boundary crossing probabilities when the original boundary~$g$ is replaced by~$f$. Upper bounds for this error have been obtained in the case when the error boundary $h = f-g$ (which can be viewed as the perturbation from~$g$) belongs to the Cameron--Martin space~\cite{Novikov1999}, and when it is Lipschitz continuous~\cite{Borovkov2005}. Beyond error bounds, the sensitivity of the non-crossing probability to perturbations of the boundary has also been studied. It was proved in~\cite{Borovkov2010} that the boundary non-crossing probability as a functional of the boundary is G\^ateaux differentiable, and an expression for the derivative was obtained. The authors in~\cite{Liang2025} used a different technique to obtain a different expression for the derivative and reconciled it with the previous expression derived in~\cite{Borovkov2010}. 

\subsection{Overview of results}

Motivated by the limited class of boundaries with explicitly known first-passage time distributions, we develop a function-space Taylor expansion around an analytically tractable reference boundary. 

Let $W := \{W_t : t \in[0,T]\}$ be a standard Brownian motion and let $g:[0,T] \to \mathbb{R}$ be a continuous boundary with $g(0) >0.$ Define the first-passage time
\begin{equation}\label{defn:tau}
	\tau_g := \inf\{t \in(0,T] : W_t \geq g(t)\}, 
\end{equation}
and let $\tau_g = \infty$ if the set is empty. The main object of interest in this paper is the first-passage time distribution function and its tail, which we denote respectively by, for $t\in [0,T],$
\begin{equation}\label{defn:F}
	F[g](t) := \Px(\tau_g \leq t), \quad \overline{F}[g](t) := \Px(\tau_g > t).
\end{equation}
For mathematical convenience, for the remainder of the paper we will express our results in terms of the distribution tail~$\overline{F}.$

In our first main result we obtain the following Taylor expansion for the first-passage time distribution tail~$\overline{F}$ in the Fr\'echet sense along the Cameron--Martin space~$\cmSpace$ (Theorem~\ref{thm:multi_linear_DF}): 
\begin{equation*}
	\overline{F}[g + h] = \overline{F}[g] + \sum_{n=1}^m\frac{1}{n!} \nabla_{(h,\ldots,h)}^n \overline{F}[g] + R_m[g,h],
\end{equation*}
with an explicit bound on the remainder $R_m[g,h],$ where $\nabla_{(h,\ldots,h)}^n$ denotes the $n$-th $\cmSpace$-derivative operator (see Section~\ref{section:taylor_series_expansion}). This extends the result in~\cite{Liang2025} to the higher-order case, where the $\cmSpace$-differentiability of the non-crossing probability $\overline{F}[g](T)$ was proved. In the special case when the boundary~$g$ is linear and the boundary perturbation~$h$ is in~$\cmSpace$, the Cameron--Martin theorem was used in~\cite{Hobson1999} to obtain a weak Taylor expansion for the infinite time horizon boundary crossing probability $F[g](\infty),$ and a different interpretation for the Taylor series coefficients was obtained. We emphasise that our results hold for general moving boundaries~$g$ and that we consider the first-passage time distribution rather than the boundary crossing probability at a fixed time point.  

Next, with the goal of obtaining an analytically tractable expression for the $\cmSpace$-derivatives of~$\overline{F},$ we obtain a representation for it in terms of the integral kernels~$\psi_t^{(n)}$ from the Wiener--It\^o chaos expansion of~$\mathbf{1}_{\tau_g > t}$ (Theorem~\ref{thm:chaos}): 
\begin{align*}
	&\nabla_{(h,\ldots,h)}^{n} \overline{F}[g](t) \nonumber \\
	&\quad = (-1)^n (n!)^2\int_{0}^{t}\int_{0}^{t_n} \cdots\int_{0}^{t_2} \psi_{t}^{(n)}(t_1,\ldots,t_n)\prod_{i=1}^n \dot{h}(t_i)\,\dif t_1\,\cdots\,\dif t_{n-1}\,\dif  t_n.
\end{align*}
This is one of the key novelties in this paper. It leads to a new interpretation of the directional derivatives of $\overline{F}[g]$ and introduces a new set of objects to the study of the moving boundary crossing problem. 

We then obtain the following relatively tractable representation for the Wiener--It\^o chaos integral kernels $\psi_t^{(n)}$ in terms of Krylov--Veretennikov--Zvonkin-type operators~\cite{Veretennikov1976,Zvonkin1975} of the form $(\frac{\partial}{\partial x}Q_t\phi)(s,x),$ where~$Q_t$ is the taboo transition operator of the Brownian motion with boundary~$g$ (Theorem~\ref{thm:psi_operator_rep}):
\begin{equation*}
	\psi_{t}^{(n)}(t_1,t_2,\ldots,t_n) = \frac{1}{n!} \biggl(Q_{t_1} \frac{\partial}{\partial x}  Q_{t_2}   \cdots   \frac{\partial}{\partial x}  Q_{t_n}   \frac{\partial}{\partial x}  Q_t 1\biggr)(0,0).
\end{equation*}
In conjunction with the above multiple integral representation for $\nabla_{(h,\ldots,h)}^n\overline{F},$ this allows us to reconcile our expression for the first $\cmSpace$-derivative of~$\overline{F}$ with the one obtained in~\cite{Liang2025} and obtain higher-order analogues. 

Finally, by combining the above results, we obtain new semi-closed-form expressions for the first two $\cmSpace$-derivatives of~$\overline{F}[g]$ in the case when~$g$ is a linear boundary (Propositions~\ref{prop:deriv_bm} and~\ref{prop:deriv_2_bm}). Based on these identities we propose a second-order Taylor approximation~\eqref{defn:F_m} for the first-passage time distribution of a general $\contSpace^2$-boundary, and quantify the error of our approximation using Theorem~\ref{thm:multi_linear_DF} (Proposition~\ref{prop:approx_m}). In the examples discussed in Section~\ref{sec:examples}, our first-order approximation is more accurate than the classical tangent approximation~\cite{Lerche1986} and its hazard-rate variant~\cite{Roberts1995} over the displayed ranges, while remaining straightforward to implement.  

\subsection{Proof outline and paper organisation}

We will now roughly sketch the steps of the proof of the main results. The starting point is the observation that, by the spatial homogeneity of Brownian motion, a boundary perturbation of the form~$h = \bm{z}^{\top}\bm{h},$ $\bm{z}\in \mathbb{R}^d,$ $\bm{h} \in \cmSpace^d$ can be transferred into the drift of the process. The drift is then `removed' by using the Cameron--Martin theorem, allowing us to express the first-passage time distribution for the perturbed boundary in terms of an expected value of a Dol\'eans exponential on the event of not crossing the original boundary. Next, we introduce a collection of auxiliary martingales through a certain family of multivariate Hermite polynomials and use them to obtain a series expansion for the Dol\'eans exponential, yielding a power series representation for $\overline{F}[g + h]$ after justifying an interchange of expectation and infinite sum. We then use the multiple Wiener--It\^o integral representation and moment bounds for the above-mentioned auxiliary martingales (Lemma~\ref{thm:mvt_HMG}) to use the coefficients of the power series above to identify the $\cmSpace$-derivatives of~$\overline{F}.$ The identities for the $\cmSpace$-derivatives and the Wiener--It\^o chaos kernels in Theorems~\ref{thm:chaos} and~\ref{thm:psi_operator_rep} are obtained using the multiple Wiener--It\^o integral representation for the auxiliary martingales and some probabilistic results for solutions to terminal-boundary-value problems for the heat equation on time-dependent domains that were obtained in~\cite{Liang2025}. 

The paper is organised as follows. In Section~\ref{section:taylor_series_expansion}, we introduce some preliminary notations regarding differentiability of Banach space functions, and the Cameron--Martin space~$\cmSpace$. In Section~\ref{sec:main_result}, we state our main results and provide the proofs in Section~\ref{sec:proofs}. We conclude with five examples in Section~\ref{sec:examples}.

\section{Preliminaries}\label{section:taylor_series_expansion}

Since we will be viewing the first-passage time distribution tail~$\overline{F}$ from~\eqref{defn:F} as a differentiable function space mapping, for the convenience of the reader we recall below some standard terminology and notations for differentiation of functions along a subspace of a Banach space. For more detail on differential calculus on Banach spaces, we refer the reader to the survey in~\cite{Averbukh1967} and the books~\cite{Dieudonne1969,Aubin1984}.  

\subsection{$\domSubSpace$-differentiability on Banach spaces}

For Banach spaces $(\domSpace, \|\cdot \|_{\domSpace})$ and $(\rangeSpace, \| \cdot\|_{\rangeSpace}),$ let $\LinSpace(\domSpace,\rangeSpace)$ denote the Banach space of continuous linear maps from $\domSpace$ to $\rangeSpace$ equipped with the operator norm. 

Suppose that~$\domSubSpace$ is a vector subspace of~$\domSpace,$ and~$\domain$ is an open subset of~$\domSpace$. We call a mapping $V : \domain \to \rangeSpace$ \textit{weakly $\domSubSpace$-differentiable at} $f \in \domain$ if there exists $L_f \in \LinSpace(\domSubSpace, \rangeSpace)$ such that, for $h \in \domSubSpace,$
\begin{equation}\label{defn:weak_h_derivative}
	V[f + \varepsilon h] - V[f] = \varepsilon L_f [h] + o(\varepsilon), \quad \varepsilon \to 0,
\end{equation}
and set 
\begin{equation*}
	\nabla_hV[f]:= L_f[h], \quad h \in \domSubSpace.
\end{equation*}
We call the mapping~$V$ \textit{strongly $\domSubSpace$-differentiable at} $f \in \domain$ if there exists $L_f\in \LinSpace(\domSubSpace, \rangeSpace)$ such that 
\begin{equation}\label{defn:strong_h_derivative}
	V[f + h] - V[f] = L_f[h] + o(\| h\|_{\domSubSpace})
\end{equation}
as $\|h\|_{\domSubSpace} \to 0.$ If $L_f$ exists in the sense of \eqref{defn:weak_h_derivative} or \eqref{defn:strong_h_derivative}, we call~$L_f$ the \textit{weak/strong $\domSubSpace$-derivative of~$V$ at~$f$}, respectively. We will use the same notation~$\nabla_hV$ for the two derivatives and clearly state the type of differentiability to distinguish them. If~$V$ is weakly/strongly $\domSubSpace$-differentiable for all $f \in \domain,$ we simply say that~$V$ is \textit{weakly/strongly $\domSubSpace$-differentiable}. For a weak $\domSubSpace$-derivative to be a strong $\domSubSpace$-derivative, it is sufficient for~$L_f$ to be continuous with respect to~$f.$ This follows from the usual proof for extending G\^ateaux differentiability to Fr\'echet differentiability (see e.g.~Lemma 1 in~\cite{Aubin1984}). 

\begin{lemma}\label{lemma:weak_to_strong_derivative}
Suppose that $V : \domain \to \rangeSpace$ is weakly $\domSubSpace$-differentiable in a neighbourhood of $f \in \domain$ and denote its weak $\domSubSpace$-derivative at~$g$ by~$L_g.$ If the mapping $g\mapsto L_g$ is continuous at~$f$ in the operator norm, then~$V$ is strongly $\domSubSpace$-differentiable at~$f,$ and its weak and strong $\domSubSpace$-derivatives coincide.
\end{lemma}

We say that the function $L : \domSubSpace^n \to \rangeSpace$ is \textit{$n$-linear and bounded} if~$L$ is linear in each argument and there exists a constant $C \geq 0$ such that
\begin{equation}\label{defn:boundedness}
	\| L[h_1,\ldots,h_n]\|_{\rangeSpace} \leq C \prod_{i=1}^n \| h_i \|_{\domSubSpace}, \quad (h_1,\ldots h_n) \in \domSubSpace^n, 
\end{equation}
and denote the space of bounded $n$-linear maps from $\domSubSpace^n$ to~$\rangeSpace$ by $\LinSpace^n(\domSubSpace , \rangeSpace).$ For $n\in \mathbb{N}_{\geq 0} ,$ we call a function~$V$ \textit{$n$-times weakly $\domSubSpace$-differentiable at $f \in \domain$} if there exists a $n$-linear bounded mapping $L_f^n \in \LinSpace^n(\domSubSpace , \rangeSpace)$ such that, for all $(h_1,\ldots,h_n) \in \domSubSpace^n,$
\begin{equation}\label{defn:nabla_hF}
	\frac{\partial^n}{\partial \varepsilon_1\cdots \partial \varepsilon_n}V \biggl[f + \sum_{i=1}^n \varepsilon_i h_i  \biggr]\Bigg|_{\varepsilon_1,\,\ldots,\varepsilon_n=0}  = L_f^n[h_1,\ldots,h_n].
\end{equation}
If~$L_f^n$ exists, we call~$L_f^n$ the \textit{$n$-th weak $\domSubSpace$-derivative} of~$V$ at~$f,$ and set 
\begin{equation*}
\nabla_{(h_1,\ldots,h_n)}^nV[f] := L_f^n[h_1,\ldots,h_n].
\end{equation*}
We define the higher-order strong $\domSubSpace$-derivatives inductively. Suppose that $V$ is $(n-1)$-times strongly $\domSubSpace$-differentiable with $(n-1)$-th strong derivative at~$f$ denoted by $L_f^{n-1} \in \LinSpace^{n-1}(\domSubSpace,\rangeSpace).$ We call a function~$V$ \textit{$n$-times strongly $\domSubSpace$-differentiable} at~$f$ if there exists $L_f^n \in \LinSpace^n(\domSubSpace , \rangeSpace)$ such that, for all $(h_1,\ldots,h_n) \in\domSubSpace^n,$
\begin{equation*}
	(L_{f + h_n}^{n-1} - L_{f}^{n - 1})[h_1,\ldots,h_{n - 1}] = L_f^n[h_1, \ldots, h_n] + r_f[ h_1, \ldots, h_n],
\end{equation*}
where
\begin{equation*}
	\sup_{\| h_1\|_{\domSubSpace} = \cdots =\| h_{n - 1}\|_{\domSubSpace}  =1}\|r_f[h_1,\ldots,h_n]\|_{\rangeSpace} = o(\| h_n\|_{\domSubSpace}),\quad \|h_n\|_{\domSubSpace} \to 0.
\end{equation*}
If $L_f^n$ exists, we call $L_f^n$ the \textit{$n$-th strong $\domSubSpace$-derivative} of~$V$ at~$f.$ If~$V[f]$ is $n$-times weakly/strongly $\domSubSpace$-differentiable for all $f \in \domain,$ we simply say that~$V$ is \textit{$n$-times weakly/strongly $\domSubSpace$-differentiable}.

\subsection{Cameron--Martin space}

For a given $T>0,$ let $\contSpace := \contSpace([0,T])$ denote the Banach space of real-valued continuous functions on $[0,T]$ endowed with the uniform norm $\| \cdot \|_{\infty}.$ 

Denote by $\cmSpace_t$ the Cameron--Martin space of absolutely continuous functions with square-integrable derivatives:
\begin{equation*}
 \cmSpace_t := \left\{h \in \contSpace([0,t]): h(s) = \mbox{$\int_0^s$}\ell(u)\,\dif u ,\, \ell\in \LpSpace^2([0,t])\right\}, \quad t \in [0,T],
\end{equation*}
and set $\cmSpace := \cmSpace_T.$ Denoting the derivative of $h$ by $\dot{h}$ (recall that $h \in \cmSpace$ are point-wise differentiable almost everywhere, see Corollary 12 in Chapter~5 in~\cite{Royden1988}), $\cmSpace_t$ is a Hilbert space when equipped with the inner product 
\begin{equation*}
	\langle h_1,h_2\rangle_{\cmSpace_t} := \int_0^t\dot{h}_1(s)\dot{h}_2(s)\,\dif s,\quad h_1, h_2\in \cmSpace_t, 
\end{equation*}  
and we set $\| h\|_{\cmSpace_t} := \langle h,h\rangle_{\cmSpace_t}^{1/2}$. 

It will be shown in the next section that~$\overline{F}$ is infinitely many times differentiable along the subspace~$\cmSpace.$ For more detail on $\cmSpace$-differentiability, see~\cite[p.~323]{IkedaWatanabe}.

\section{Main results}\label{sec:main_result}

The goal of this paper is to obtain a new series expansion for the first-passage time distribution tail~$\overline{F}$ defined in~\eqref{defn:F} for the Brownian motion with a time-dependent boundary. We start by showing in Theorem~\ref{thm:multi_linear_DF} that the distribution tail~$\overline{F}$ admits a Taylor series representation. As an intermediate step towards a more explicit expression for the $\cmSpace$-derivatives of~$\overline{F}$, we establish in Theorem~\ref{thm:chaos} the Wiener--It\^o chaos expansion for $\mathbf{1}_{\tau_g > t},$ and an identity for its kernel functions (which we will refer to as the \emph{chaos kernels} for brevity) in terms of the $\cmSpace$-derivatives of~$\overline{F}.$ In Theorem~\ref{thm:psi_operator_rep} we obtain an explicit representation for the chaos kernels in terms of Krylov--Veretennikov--Zvonkin-type operators, which are relatively tractable in the case when the taboo transition density of the Brownian motion is known in closed-form, e.g.\ in the linear boundary case. Using this representation, we propose an approximation for~$\overline{F}$ based on semi-closed-form expressions for its first two $\cmSpace$-derivatives obtained in Propositions~\ref{prop:deriv_bm} and~\ref{prop:deriv_2_bm}, and quantify the approximation error in Proposition~\ref{prop:approx_m}. Proofs of the results obtained in this section will be presented in Section~\ref{sec:proofs}.

For the remainder of the paper we will assume that the boundary~$g$ starts above the Brownian motion, i.e.\ is in the following open subset of~$\contSpace$:
\begin{equation*}
\contSpaceZero := \{ f \in \contSpace  : f(0) > 0\}.
\end{equation*}
The following continuity property formalises the boundary-to-distribution mapping viewpoint taken in this paper. 

\begin{lemma}\label{lemma:F_cont}
The function $\overline{F}$ defined in~\eqref{defn:F} is a continuous map from~$ \contSpaceZero $ to~$ \contSpace .$
\end{lemma}

\subsection{Taylor expansions and Wiener--It\^o chaos expansions}\label{sec:taylor_series}

The first main result of the paper is a Taylor expansion for the first-passage time distribution tail~$\overline{F},$ with an explicit bound for the truncation error of the expansion. 
 
\begin{theorem}\label{thm:multi_linear_DF}
The mapping~$\overline{F}: \contSpaceZero \to \contSpace $ defined in~\eqref{defn:F} is infinitely many times strongly $\cmSpace$-differentiable. Furthermore, for all $m\in \mathbb{N},$ $g \in \contSpaceZero,$ $h \in \cmSpace,$ one has
\begin{equation}\label{eqn:taylor_series}
	\overline{F}[g+h] = \overline{F}[g] + \sum_{n=1}^{m}\frac{1}{n!}\nabla_{(h,\ldots, h)}^n\overline{F}[g] +R_m[g,h],
	\end{equation}
where the remainder $R_m : \contSpaceZero \times \cmSpace \to \contSpace $ has the following upper bound\textnormal{:}
\begin{equation}\label{eqn:error_remainder_taylor}
	\lvert R_m[g,h](t)\rvert \leq \overline{F}[g](t)^{1/2}\sum_{n=m+1}^{\infty}\frac{\| h\|_{\cmSpace_t}^n}{\sqrt{n!}} , \quad t \in [0,T].
\end{equation}
\end{theorem} 

\begin{remark}
Aside from being a method of obtaining higher-order approximations for $\overline{F}[g+h],$ a useful feature of the Taylor expansion~\eqref{eqn:taylor_series} of~$\overline{F}$ is that the series summands are proportional to the $\cmSpace$-derivatives of~$\overline{F}$ and hence have direct interpretations as sensitivities. These derivatives arise naturally in several contexts, including the local asymptotic power of the Kolmogorov--Smirnov test~\cite{Andel1967,Gur2019}, corrections to diffusion approximations in sequential analysis~\cite{Hogan1984}, and refinements of numerical methods for computing boundary crossing probabilities~\cite{Borovkov2010, Liang2023, Liang2025}. 

A related Taylor expansion for $\overline{F}[g+h](\infty)$ was obtained in~\cite{Hobson1999} in the special case when~$g$ is linear and $h \in \cmSpace_{\infty}.$ However, its coefficients were not identified as $\cmSpace$-derivatives of~$\overline{F},$ and no bound for the truncation error was provided.
\end{remark}
 
The next theorem establishes the Wiener--It\^o chaos expansion for $\mathbf{1}_{\tau_g > t}$ and relates its chaos kernels with the $\cmSpace$-derivatives of~$\overline{F}.$ For an introduction to Wiener--It\^o chaos expansions, we refer the reader to~\cite[Ch.~1]{Nualart2006}. 

For $n\in\mathbb{N},$ denote the ordered $n$-simplices by
\begin{equation*}
	\Delta_n(t) := \{(t_1,\ldots, t_n)  : 0 < t_1 < t_2 < \cdots < t_n <t \}, \quad t \in (0,T].
\end{equation*}
For $f \in \LpSpace^2(\Delta_n(t)),$ we denote the $n$-multiple Wiener--It\^o integral of~$f$ by
\begin{equation*}
	I_{t}^{(n)}(f) := n!\int_0^{t}\int_0^{t_n}\cdots\int_0^{t_2} f(t_1,\ldots,t_n)\,\dif W_{t_1}\cdots \,\dif W_{t_{n-1}}\,\dif W_{t_n}, \quad t \in [0,T],
\end{equation*}
where $I_t^{(0)}$ is the identity mapping on the constants and~$W$ is the same Brownian motion used in the definition of the first-passage time~\eqref{defn:tau}.  

\begin{theorem}\label{thm:chaos}
For $g \in \contSpaceZero$ and $t\in(0,T],$ there exists a unique sequence of functions $\psi_{t}^{(n)}\in \LpSpace^2(\Delta_n(t)),$ $n \in \mathbb{N},$ such that
\begin{equation}\label{eqn:series_indicator}
	\mathbf{1}_{\tau_g > t}  = \sum_{n=0}^{\infty} I_{t}^{(n)}(\psi_{t}^{(n)}),
\end{equation}
where $\psi_t^{(0)}= \overline{F}[g](t)$ and the series on the right-hand side of~\eqref{eqn:series_indicator} converges in the mean-square sense. Letting $\kappa_t(s) := s \wedge t,$ $s\in [0,T],$ one has, for $n\in\mathbb{N},$
\begin{equation}\label{eqn:psi_n_fla}
		\psi_{t}^{(n)}(t_1,\ldots,t_n)	  =   \frac{(-1)^n}{n!} \frac{\partial^n}{\partial t_1\cdots\partial t_n} \nabla_{(\kappa_{t_1},\ldots,\kappa_{t_n})}^n \overline{F}[g](t) 
\end{equation}
a.e.\ on $(t_1,\ldots,t_n)\in \Delta_n(t).$ Moreover,  
\begin{equation}\label{eqn:dF_multiple_integral}
\nabla_{(h,\ldots,h)}^{n} \overline{F}[g](t)  = (-1)^n (n!)^2 \int_{\Delta_n(t)}\psi_t^{(n)}(\bm{t})\dot{h}^{\otimes n}(\bm{t})\,\dif \bm{t},
\end{equation}
where $\dot{h}^{\otimes n}(t_1,\ldots,t_n):= \prod_{i=1}^n \dot{h}(t_i).$
\end{theorem}

\begin{remark}
For sufficiently regular Brownian functionals, one can obtain expressions for the chaos kernels~$\psi_t^{(n)}$ in terms of iterated Malliavin derivatives using the Taylor--Stroock formula~\cite{Stroock1987}. This approach does not apply directly here since $\mathbf{1}_{\tau_g > t}$ is not Malliavin differentiable, a fact which can be verified by using an elementary result by Sekiguchi and Shiota~\cite{Sekiguchi1985} (see~\cite[Prop.~1.2.6]{Nualart2006} for a more accessible reference). Theorem~\ref{thm:chaos} instead recovers the chaos kernels from the $\cmSpace$-derivatives of the map~$g \mapsto \overline{F}[g].$
\end{remark}

\begin{remark}
Our representation~\eqref{eqn:psi_n_fla} for the chaos kernels can be viewed as a higher-order extension of the Riesz–Fr\'echet representation for the G\^ateaux derivative of~$\overline{F}$ (see Remark 7 in~\cite{Liang2025}). Moreover, approximating the mixed derivatives in~\eqref{eqn:psi_n_fla} using finite differences yields numerical approximations of the chaos kernels. We demonstrate this procedure in Example~\ref{ex:fpt_to_chaos} in Section~\ref{sec:examples}.  
\end{remark}

In the next theorem we obtain a representation for the chaos kernels~$\psi_{t}^{(n)}$ which is more suitable for efficient numerical computation when one has access to the taboo transition density. This will allow~\eqref{eqn:dF_multiple_integral} to become a useful formula for computing the $\cmSpace$-derivatives of~$\overline{F}$.

For $(s,x) \in [0,T] \times \mathbb{R},$ let
\begin{equation*}\label{diffusion_process_Y}
	W_t^{s,x} := x + W_{t-s}, \quad t \in [s,T].
\end{equation*} 
Since our proof will use technical results from parabolic partial differential equations, we will assume that~$g$ is twice-continuously differentiable, i.e.\ $g \in  \contSpace^{2} \equiv \contSpace^{2}([0,T]).$ Define the family of time-space domains
\begin{equation*}
	D_t := \{(s,x) : s \in (0,t) ,\ x < g(s)\} \subset \mathbb{R}^2 , \quad t\in (0,T],
\end{equation*}
and the first-passage times
\begin{equation*}
	\tau_g^{s,x} := \inf\{u > s: W_u^{s,x}\geq g(u)\}, \quad (s,x) \in \text{cl}(D_T),
\end{equation*}
where $\text{cl}(\,\cdot\,)$ denotes the set closure. For bounded functions $\phi : (0,\infty) \times \mathbb{R}\to \mathbb{R},$ we denote the taboo transition operator~$Q_t$ of the Brownian motion~$W$ by
\begin{equation*}\label{defn:U}
	(Q_{t}\phi)(s,x) := \Ex\, \phi(t,W_t^{s,x})\bm{1}_{\tau_g^{s,x}> t}, \quad (s,x) \in D_t.
\end{equation*}
For brevity, we will write $(Q_s\Dx Q_t \phi)(r,x):= (Q_s(\Dx(Q_t\phi)))(r,x).$ 

\begin{theorem}\label{thm:psi_operator_rep}
Let $g \in \contSpaceZero \cap \contSpace^2$ and $t\in(0,T].$ One has, for $n \geq 1,$  
\begin{equation}\label{eqn:equiv_psi}
	\psi_{t}^{(n)}(t_1,t_2,\ldots,t_n) = \frac{1}{n!} \biggl(Q_{t_1} \frac{\partial}{\partial x}  Q_{t_2}   \cdots   \frac{\partial}{\partial x}  Q_{t_n}   \frac{\partial}{\partial x}  Q_t 1\biggr)(0,0)
\end{equation}
a.e.\ on $(t_1,\ldots,t_n) \in \Delta_n(t).$
\end{theorem}

\begin{remark}
Representation~\eqref{eqn:equiv_psi} allows one to reconcile the first chaos kernel from Theorem~\ref{thm:chaos} with the kernel function~$\psi$ appearing in the first G\^ateaux derivative of~$\overline{F}$ obtained in~\cite{Liang2025}. Expressions of the form on the right-hand side of~\eqref{eqn:equiv_psi} with the taboo transition operators~$Q_t$ replaced by the `free' transition operators~$P_t,$ where $(P_t\phi)(s,x) := \Ex\,\phi(t,W_t^{s,x}),$ arise in the Krylov--Veretennikov--Zvonkin~\cite{Veretennikov1976,Zvonkin1975} representation for strong solutions of stochastic differential equations. An analogous representation to~\eqref{eqn:equiv_psi} in the case of a two-sided constant boundary was obtained in the study of the Onsager--Machlup functional~\cite{Carmona1992}. 
\end{remark}

\begin{remark}
It was proved in Lemma 5 in~\cite{Liang2025} that if~$g \in \contSpaceZero\cap \contSpace^2,$ the taboo transition density
\begin{equation*}
	q(s,x;t,y) := \Px(W_t^{s,x}\in \dif y,\tau_g^{s,x}> t)/\dif y
\end{equation*} 
exists for $0 \leq s < t \leq T,$ $x < g(s),$ $y <g(t).$ Therefore, one can express the taboo transition operator as follows:
\begin{equation}\label{eqn:operator_taboo_transition}
	(Q_t\phi)(s,x) = \int_{-\infty}^{g(t)}q(s,x;t,y)\phi(t,y)\,\dif y,\quad (s,x) \in D_t,\ t \in (0,T].
\end{equation}
This allows one to view the right-hand side of~\eqref{eqn:equiv_psi} as an $n$-dimensional integral. 
\end{remark}

\begin{remark} 
By repeatedly applying integration by parts, one can obtain from~\eqref{eqn:equiv_psi}
\begin{equation}\label{eqn:taboo_psi}
\begin{aligned}
	\psi_{t}^{(n)}(t_1,t_2,\ldots,t_n)&= \frac{(-1)^n}{n!}\int_{-\infty}^{g(t_1)} \cdots  \int_{-\infty}^{g(t_{n-1})}\int_{-\infty}^{g(t_n)} (Q_t1)(t_n,x_n) \\
		&\quad \times \prod_{i=1}^n \frac{\partial}{\partial x_i}q(t_{i-1},x_{i-1};t_i, x_i)\,\dif x_n\,\dif x_{n-1}\cdots \,\dif x_1,
\end{aligned}
\end{equation}
a.e.\ on $(t_1,\ldots,t_n) \in \Delta_n(t).$ Expressions of the form on the right-hand side of~\eqref{eqn:taboo_psi} with the taboo transition densities~$q$ replaced by the `free' transition densities of the Brownian motion can be traced back to Isobe and Sato~\cite{Isobe1983}, see also~\cite{Hu1997}. 
\end{remark}

\subsection{Second-order approximation of the first-passage time distribution}\label{sec:fpt_approx}

In this section we propose a new approximation for the first-passage time distribution tail~$\overline{F}$ based on the Taylor expansion~\eqref{eqn:taylor_series} from Theorem~\ref{thm:multi_linear_DF} and derive semi-closed-form expressions for the $\cmSpace$-derivatives of~$\overline{F}$ using~\eqref{eqn:dF_multiple_integral} and~\eqref{eqn:equiv_psi} from Theorems~\ref{thm:chaos} and~\ref{thm:psi_operator_rep}, respectively. 

For $t \in [0,T],$ we approximate $g \in \contSpaceZero \cap \contSpace^2$ with a family of linear boundaries $\ell_t : [0,t] \to \mathbb{R}$ with matching endpoints with~$g$ on the intervals $[0,t],$ i.e.
\begin{equation}\label{defn:g_t}
	\ell_t(s) := a + b_ts , \quad 0 \leq s \leq t \leq T,
\end{equation}
where
\begin{equation*}
a := g(0),\quad  b_t := \begin{cases}
					\frac{g(t) -g(0)}{t} , & t>0,\\
					g'(0), & t =0.
				\end{cases}
\end{equation*}
We introduce the following three levels of approximations for the first-passage time distribution tail~$\overline{F}[g]$: for $t \in (0,T],$ 
\begin{equation}\label{defn:F_m}
	\overline{F}_m(t)  := \begin{cases}
		\overline{F}[\ell_t](t) , &  m =0, \\
		\overline{F}[\ell_t](t)  + \nabla_{h_t}\overline{F}[\ell_t](t)  ,& m = 1,\\
		\overline{F}[\ell_t](t)  + \nabla_{h_t}\overline{F}[\ell_t](t) + \mbox{$\frac{1}{2}$}\nabla_{(h_t,h_t)}^2\overline{F}[\ell_t](t)  ,& m = 2,
	\end{cases}
\end{equation}
where the $\cmSpace$-derivatives are computed using Propositions~\ref{prop:deriv_bm} and~\ref{prop:deriv_2_bm} below. The error arising from this approximation is quantified in Proposition~\ref{prop:approx_m} below.  

\begin{remark}
The reason why we chose the linear function~$\ell_t$ in~\eqref{defn:g_t} is because the average slope~$b_t$ is simple to calculate and it minimises the $\cmSpace$-norm of the boundary approximation error in the class of linear boundaries with matching values at $t=0.$ Optimal slopes other than~$b_t$ can be derived under different optimality criteria for the boundary approximation error, but they tend to be more complicated and thus we leave it for future research. 
\end{remark}

Denote the density of the first-passage time by $f_{\tau_g}(t) := \frac{\dif}{\dif t} F[g](t),$ $t\in (0,T),$ whose existence is guaranteed since $g \in \contSpace^1((0,T)),$ see e.g.\ Theorem~4.1 in \cite{Peskir2002}. For $t >0,$ $x \in \mathbb{R},$ set
\begin{equation}\label{defn:varphi}
	\varphi_t(x) := \frac{1}{\sqrt{2\pi t}}\exp\biggl\{-\frac{x^2}{2t}\biggr\},\quad \Phi_t(x):=  \int_{-\infty}^x \varphi_t(z)\,\dif z.
\end{equation}
 
\begin{proposition}\label{prop:deriv_bm}
For $a > 0,$ $b\in \mathbb{R},$ let $g(s) := a+bs,$ $s\geq 0,$ and $h \in \cmSpace.$ Then  
\begin{equation}\label{eqn:dF_1_BM}
\nabla_{h}\overline{F}[g](t) = \int_0^t h(s) \mbox{$\frac{\partial}{\partial s}$}\psi_{t}^{(1)}(s)\,\dif s 
\end{equation}
for $t \in (0,T],$ where, for a.e.\ $s \in(0,t)$, one has
\begin{equation}\label{defn:first_chaos_linear}
\psi_{t}^{(1)}(s)	= \int_s^t (\Dx Q_t1)(u,g(u)-)f_{\tau_g} (u)\,\dif u ,
\end{equation}
with $f_{\tau_g}(s) = \mbox{$\frac{a}{s}$}\varphi_s(a + bs)$ and 
\begin{equation}\label{eqn:bdry_grad_bm}
	  (\Dx Q_t1) (s,g(s)-) = -2\varphi_{t-s}(b(t-s)) - 2b \Phi_{t-s}(b(t-s)).
\end{equation}
\end{proposition}

\begin{proposition}\label{prop:deriv_2_bm}
For $a > 0,$ $b\in \mathbb{R},$ let $g(s) := a+bs,$ $s\geq 0,$ and $h \in \cmSpace.$ Then   
\begin{equation}\label{eqn:dF_2_BM}
\begin{aligned}
	\nabla_{(h,h)}^2 \overline{F}[g](t) &= -4\int_0^t h(t_2)\dot{h}(t_2) \psi_{t}^{(2)}(t_2-,t_2)\,\dif t_2 \\
	&\quad -4\int_0^{t}\int_0^{t_2}h(t_2)\dot{h}(t_1) \mbox{$\frac{\partial}{\partial t_2}$}\psi_{t}^{(2)}(t_1,t_2)\,\dif t_1\,\dif t_2
	\end{aligned}
\end{equation}
for $t\in (0,T],$ where, for a.e.\ $0 < t_1<  t_2 < t,$ one has
\begin{equation}\label{defn:second_chaos_linear}
	\psi_{t}^{(2)}(t_1,t_2)   =  -\int_{t_2}^t f (t_1, s)\mbox{$\frac{\partial}{\partial s}$}\psi_t^{(1)}(s)\,\dif s,
\end{equation}
with $\psi_t^{(1)}$ given by~\eqref{defn:first_chaos_linear} and
\begin{equation}\label{eqn:K(ru)}
\begin{aligned}
	f(r,s) &:= \frac{b}{2} + \frac{1}{2} \bigg(\frac{a}{s} - \frac{1}{a} \bigg)(2\Phi_c(a) -1)+  \frac{c}{s}\varphi_c(a)  ,
\end{aligned}
\end{equation}
for $0< r < s,$ where $c := rs/(s-r).$
\end{proposition}

\begin{remark}
From~\eqref{eqn:dF_multiple_integral} and~\eqref{eqn:taboo_psi}, one can see that the $n$-th $\cmSpace$-derivative $\nabla_{(h,\ldots,h)}^n \overline{F}[g]$ can be expressed as a $(2n + 1)$-dimensional integral. Propositions~\ref{prop:deriv_bm} and~\ref{prop:deriv_2_bm} show that one can reduce the dimensionality of the integrals needed to compute the $\cmSpace$-derivatives by using integration by parts. Proposition~\ref{prop:deriv_bm} reduces the dimensionality of the integral required to compute $\nabla_h \overline{F}[g]$ from three to one. Proposition~\ref{prop:deriv_2_bm} reduces the dimensionality of the integral required to compute $\nabla_{(h,h)}^2 \overline{F}[g]$ from five to two since $\frac{\partial}{\partial t_2}\psi_t^{(2)}(t_1,t_2)$ has a closed-form expression due to~\eqref{defn:second_chaos_linear}.
\end{remark}

\begin{remark}
The term $(\Dx Q_t1)(s,g(s)-)$ in the integral expression for the first chaos kernel~\eqref{defn:first_chaos_linear} appears in the context of boundary overshoot corrections for boundary crossing probabilities for random walk approximations of diffusion processes~\cite{Gobet2010,Keener2013}, and time-dependent barrier option pricing~\cite{Mijatovic2010}. In the process of reconciling the different expressions for $\nabla_h\overline{F}[g](T)$ from~\cite{Borovkov2010} and~\cite{Liang2025}, it was proved in~\cite{Liang2025} that $(\Dx Q_t1)(s,g(s)-)$ can be expressed in terms of the expected value of a functional of a Brownian meander (see Remark 2 and Lemma~4 in~\cite{Liang2025}). Related Bessel-bridge expressions appear in the context of integration by parts formulae for restricted Wiener measures and their close relatives~\cite{Funaki2007,Zambotti2002}. 
\end{remark}

In the next proposition we use Theorem~\ref{thm:multi_linear_DF} to obtain error bounds arising from using approximation~\eqref{defn:F_m}, which are useful for both the~$T\to 0$ and $T \to \infty$ cases.  Denote the error boundary by
\begin{equation}\label{defn:h_t}
	h_t(s) := g(s) - \ell_t(s), \quad 0 \leq s \leq t \leq T.
\end{equation}
 
\begin{proposition}\label{prop:approx_m}
Suppose that $g \in \contSpaceZero \cap \contSpace^2$ and recall the definition of~$\overline{F}_m$ and~$h_t$ from~\eqref{defn:F_m} and~\eqref{defn:h_t}, respectively. For $m =0,1,2,$ one has, for $t \in(0,T],$
	\begin{equation}\label{ineq:bound_Fm}
	\lvert \overline{F}[g](t) - \overline{F}_m(t) \rvert \leq \overline{F}[\ell_t](t)^{1/2}\sum_{n= m+1}^{\infty} \frac{\| h_t\|_{\cmSpace_t}^n}{\sqrt{n!}},
	\end{equation}
where $ \| h_t\|_{\cmSpace_t} \leq \pi^{-1} t \| \ddot{g}\|_{\LpSpace^2((0,t))} <\infty.$ If the boundary~$g$ admits the representation $g(s) = a+bs +h(s),$ $s\geq 0,$ for some $a>0,$ $b \in \mathbb{R}$ and $h \in \cmSpace_{\infty},$ then 
\begin{equation*}
	 \sup_{t >0}\| h_t\|_{\cmSpace_t} \leq \| h\|_{\cmSpace_{\infty}}.
\end{equation*}
\end{proposition}

In Section~\ref{sec:examples}, we will compare the accuracy of our proposed approximation~\eqref{defn:F_m} with the classical tangent approximation and its hazard-rate variant in two examples. 

\section{Proofs}\label{sec:proofs}

In this section we provide proofs for the claims made in Section~\ref{sec:main_result}. 

\begin{proof}[Proof of Lemma~\ref{lemma:F_cont}]
For $g \in \contSpaceZero,$ set
\begin{equation}\label{defn:B_g}
B_g(t) := \{w \in \contSpace : w(s) < g(s),\ s\in[0,t]\}, \quad t\in[0,T].
\end{equation}
Clearly, 
\begin{equation*}
	\overline{F}[g](t) \equiv \Px(\tau_g > t) = \Px(W \in B_g(t)), \quad t \in [0,T].
\end{equation*}
For $\delta \in (0,g(0)),$ let~$f \in \contSpaceZero$ satisfy $\| g - f \|_{\infty} \leq \delta.$ Since $g(t) -\delta \leq f(t) \leq g(t) + \delta,$ it follows that $B_g(t)\Delta B_f(t) \subseteq B_{g+\delta}(t) \backslash B_{g-\delta}(t)$ for all $t\in[0,T]$, where~$\Delta$ denotes the symmetric set difference. Hence, for all $t \in [0,T],$
\begin{align*}
	\lvert \Px(W \in B_g(t)) -  \Px(W \in B_f(t))\rvert&\leq  \Px(W \in B_g(t)\Delta B_{f}(t)  ) \leq  \phi_{\delta}(t),
\end{align*}
where 
\begin{equation*}
	\phi_{\delta}(t) := \Px(W\in B_{g+\delta}(t)\backslash B_{g-\delta}(t))  = \Px(\tau_{g+\delta} > t ) - \Px(\tau_{g-\delta} > t), \quad t \in[0,T].
\end{equation*}
For all $\delta \in (0,g(0)),$ $\phi_{\delta}(t)$ is continuous in $t \in [0,T]$ since~$\Px (\tau_{g\pm\delta} =t) =0$ because~$g$ is continuous and $g(0)\pm \delta >0.$ Since the sets $B_{g+\delta}(t)\backslash B_{g-\delta}(t)$ decrease as $\delta\downarrow 0,$ we have by continuity from above
\begin{equation*}
	\lim_{\delta \downarrow 0}\phi_{\delta}(t) = \Px(\cap_{\delta >0} (B_{g+\delta}(t) \backslash B_{g-\delta}(t))), \quad t \in [0,T].
\end{equation*}
Since~$g$ is continuous it is not hard to see that the event on the right-hand side, i.e.\ the set of paths that touch the boundary~$g$ but never cross it before time $t,$ has probability zero (see e.g.\ Proposition 3.4 in~\cite{Antunovic2011}). Choosing a sequence $\delta_k\downarrow 0,$ the sequence $\{\phi_{\delta_k}\}$ is decreasing, each element being continuous, and converges pointwise to zero. By Dini's theorem (see e.g.~\cite[Theorem 7.13]{Rudin1976}) it follows that $\sup_{t \in [0,T]}\lvert \phi_{\delta_k}(t)\rvert \to 0,$ and thus the mapping $g \mapsto \Px(W \in B_g(\,\cdot\,))$ from $\contSpaceZero$ to $\contSpace$ is continuous with respect to the uniform norm.
\end{proof}

\subsection{Taylor expansions and Wiener--It\^o chaos expansions}

To prove Theorems~\ref{thm:multi_linear_DF}--\ref{thm:psi_operator_rep} we will first introduce a family of auxiliary martingales which are extensions of the univariate Hermite polynomials applied to the Brownian motion, (for the latter, see~\cite{Segall1976}, or Proposition~3.8 in~\cite[p.~151]{RevuzYor} for a more accessible reference). 

For $d\in \mathbb{N}_{>0}$ and a multi-index $\bm{\alpha} := (\alpha_1,\ldots,\alpha_d) \in \mathbb{N}_{\geq 0}^d$ (all the vectors in what follows are column vectors),  let 
\begin{equation*}
	\lvert \bm{\alpha}\rvert := \alpha_1 + \alpha_2 +\cdots +\alpha_d, \quad \bm{\alpha}! := \alpha_1!\alpha_2!\cdots \alpha_d!,
\end{equation*}
and for a vector $\bm{z} := (z_1,\ldots,z_d)\in\mathbb{R}^d$ we set
\begin{equation*}
	\bm{z}^{\bm{\alpha}} := z_1^{\alpha_1}z_2^{\alpha_2}\cdots z_d^{\alpha_d}.
\end{equation*}

For a covariance (symmetric and positive semidefinite) matrix $\boldsymbol{\Sigma} \in \mathbb{R}^{d\times d}$ with~$\Sigma_{i,j}$ denoting its $(i,j)$-entry, define the multivariate Hermite polynomials $\dHe_{\bm{\alpha}}(\bm{x} ;\bm{\Sigma}),$ $\bm{x} \in \mathbb{R}^d,$ $\bm{\alpha}\in \mathbb{N}_{\geq 0}^d,$ associated with the covariance matrix~$\bm{\Sigma}$ through the  following generating function (cf.\ \cite[(A.6)]{Takemura1988}): for $\bm{z} \in \mathbb{R}^d,$
\begin{equation}\label{gen_func}	 
	\Lambda(\bm{z},\bm{x};\bm{\Sigma}) :=\exp\bigl\{\bm{z}^{\top} \bm{x} - \mbox{$\frac{1}{2} $}\bm{z}^{\top} \bm{\Sigma} \bm{z}\bigr\}= \sum_{\bm{\alpha} \in \mathbb{N}_{\geq 0}^d}\frac{\bm{z}^{\bm{\alpha}}}{\bm{\alpha}!}\dHe_{\bm{\alpha}}(\bm{x} ; \bm{\Sigma}) .
\end{equation}
In the special case when~$\bm{\alpha}$ is equal to
\begin{equation}\label{defn:bm_1}
	 \mathbf{1}_d := (1,\ldots,1)\in \mathbb{N}^d,
\end{equation}
one has (see \cite[(A.16)]{Takemura1988})
\begin{align*}
	\dHe_{\mathbf{1}_1}(\bm{x} ; \bm{\Sigma}) &= x_1,\\
	\dHe_{\mathbf{1}_2}(\bm{x} ; \bm{\Sigma}) &= x_1x_2 - \Sigma_{1,2},\\
	\dHe_{\mathbf{1}_3}(\bm{x} ; \bm{\Sigma}) &= x_1x_2x_3 -x_1\Sigma_{2,3} -x_2\Sigma_{1,3} - x_3 \Sigma_{1,2}.
\end{align*}

Let $\bm{h} := (h_1,\ldots,h_d) \in \cmSpace^d.$ For $i=1,\ldots,d,$ define 
\begin{equation}\label{defn:M^h}
	M_t^{(i)} := -\int_0^t \dot{h}_i(s)\,\dif W_s, \quad t\in[0,T],
\end{equation}
which has quadratic variation $[M^{(i)}]_t = \| h_i\|_{\cmSpace_t}^2.$ Define the vector-valued martingale 
\begin{equation}\label{defn:vector_mg}
	\bm{M}_t^{\bm{h}} := (M_t^{(1)},\ldots,M_t^{(d)}),\quad t \in[0,T],
\end{equation}
and denote its quadratic covariation matrix by $\bm{C}_t^{\bm{h}} := (C_t^{(i,j)}),$ where its elements are given by
\begin{equation*}
C_t^{(i,j)} = [M^{(i)}, M^{(j)}]_t= \langle h_i,h_j\rangle_{\cmSpace_t},\quad i,j \in \{1,\ldots,d\}.
\end{equation*}

\begin{remark}\label{remark:compare_malliavin}
In the usual Hermite--Fourier construction of Wiener chaos expansions, one chooses an orthonormal basis $\{e_i : i \geq 1\}$ of $\LpSpace^2([0,T])$ and expands in products of univariate Hermite polynomials evaluated at the independent Gaussian random variables $\xi_i := \int_0^T e_i(s)\,\dif W_s$ (see e.g.\ Proposition 1.1.1 in~\cite{Nualart2006}). In this paper, the Gaussian vector $\bm{M}_t^{\bm{h}}$ arises directly from the perturbation~$\bm{h}$ and its components do not need to be independent. The polynomials~$\dHe_{\bm{\alpha}}$ allow us to retain these coordinates and their covariance matrix without first orthonormalising the perturbation directions. When~$\bm{\Sigma}$ is the identity matrix, they reduce to the usual product form.
\end{remark}

\begin{lemma}\label{thm:mvt_HMG}
The process $L^{\bm{\alpha},\bm{h}} := \{L_t^{\bm{\alpha},\bm{h}} : t\in [0,T]\},$ $\bm{\alpha} \in \mathbb{N}_{\geq 0}^d,$ $\bm{h} \in \cmSpace^d,$ defined by
\begin{equation}\label{defn:L_a_h}
	L_t^{\bm{\alpha},\bm{h}} := \dHe_{\bm{\alpha}}(\bm{M}_t^{\bm{h}}; \bm{C}_t^{\bm{h}}), \quad t\in [0,T], 
\end{equation}
is a continuous martingale and has the following representation\textnormal{:}
\begin{equation}\label{multiple_integral_HMG}
L_t^{\bm{\alpha},\bm{h}} = \sum_{\bm{\sigma} \in \mathcal{P}(\bm{\alpha})} \bm{\alpha}! \int_0^{t}\int_0^{t_{\lvert \bm{\alpha}\rvert }}\cdots \int_0^{t_2}  \,\dif M_{t_1}^{(\sigma_1)}\cdots \,\dif M_{t_{\lvert \bm{\alpha}\rvert-1}}^{(\sigma_{\lvert\bm{\alpha}\rvert-1})}\,\dif M_{t_{\lvert \bm{\alpha}\rvert}}^{(\sigma_{\lvert\bm{\alpha}\rvert})}, \quad t\in [0,T],
\end{equation}
where $\bm{\sigma} := (\sigma_1, \sigma_2, \ldots,\sigma_{\lvert \bm{\alpha}\rvert}) \in \mathbb{N}^{\lvert\bm{\alpha} \rvert}$ and $\mathcal{P}(\bm{\alpha})$ is defined as the set of unique permutations of $\{1,\ldots,d\}$ with respective multiplicities given by~$\bm{\alpha}$. Moreover,
\begin{equation}\label{square_integrable_HMG}
	\sup_{t \in [0,T]}\Ex(L_t^{\bm{\alpha},\bm{h}} )^2 \leq \lvert \bm{\alpha} \rvert! \prod_{k=1}^{d}\| h_k\|_{\cmSpace}^{2 \alpha_k}.
\end{equation}
\end{lemma}

\begin{remark}
In the Gaussian Hilbert space literature (see e.g.~\cite{Janson1997}), the random variable $\Lambda(\bm{z},\bm{M}_t^{\bm{h}}; \bm{C}_t^{\bm{h}})$ is called the \emph{Wick exponential of $ \bm{z}^{\top}\bm{M}_t^{\bm{h}} $}, while the random variables $\{L_t^{\bm{\alpha},\bm{h}}:\bm{\alpha} \in \mathbb{N}_{\geq 0}^{d}\}$ are called \emph{Wick monomials}.
\end{remark}

\begin{proof}[Proof of Lemma~\ref{thm:mvt_HMG}]
It can be verified from~\eqref{gen_func} that
\begin{equation}\label{eqn:gen_func_PDE}
\begin{aligned}
	( \mbox{$ \frac{\partial^2}{\partial x_i \partial x_j} + \frac{\partial}{\partial  \Sigma_{i,j}} $})  \Lambda(\bm{z},\bm{x};\bm{\Sigma}) = 0,& \quad i=1,\ldots,d,\ i < j,\\
	( \mbox{$ \frac{1}{2}\frac{\partial^2}{\partial x_i^2} + \frac{\partial}{\partial \Sigma_{i,i}} $} )  \Lambda(\bm{z},\bm{x};\bm{\Sigma}) =0, & \quad i=1,\ldots,d.
\end{aligned}
\end{equation}
Taking partial derivatives $\frac{\partial^{\lvert \bm{\alpha}\rvert }}{\partial z_1^{\alpha_1} \cdots \partial z_d^{\alpha_d}}$ on both sides of~\eqref{eqn:gen_func_PDE} and interchanging the order of derivatives and evaluating at $\bm{z} = \bm{0},$ it follows that, for $\bm{\alpha} \in \mathbb{N}_{\geq 0}^d,$
\begin{equation}\label{covar_pde}
\begin{aligned}
		(\mbox{$\frac{\partial^2}{\partial x_i \partial x_j} + \frac{\partial}{\partial  \Sigma_{i,j}}$}) \dHe_{\bm{\alpha}}(\bm{x};\bm{\Sigma}) &=0,\quad i=1,\ldots,d,\ i< j, \\
(\mbox{$\frac{1}{2}\frac{\partial^2}{\partial x_i^2} + \frac{\partial}{\partial \Sigma_{i,i}}$}) \dHe_{\bm{\alpha}}(\bm{x};\bm{\Sigma}) &=0,\quad i=1,\ldots,d.
\end{aligned}
\end{equation}

Using It\^o's formula for the process $\dHe_{\bm{\alpha}}(\bm{M}_t^{\bm{h}};\bm{C}_t^{\bm{h}} )$ and the following recurrence relation (see~\cite[(A.9)]{Takemura1988}):
\begin{equation*}
	\mbox{$\frac{\partial}{\partial x_i}$}\dHe_{\bm{\alpha}}(\bm{x} ; \bm{\Sigma}) =\alpha_i \dHe_{\bm{\alpha}-\bm{e}_{i}}(\bm{x} ; \bm{\Sigma}), \quad i =1,\ldots,d,
\end{equation*}
where $\bm{e}_i$ denotes the column vector of length~$d$ containing zeros except for a 1 in the $i$-th entry, one has, for $t\in [0,T],$
\begin{align}\label{recursion_mvt_HMG}
L_t^{\bm{\alpha},\bm{h}} &= \sum_{i=1}^d\int_0^t \mbox{$\frac{\partial}{\partial x_i}$}\dHe_{\bm{\alpha}}(\bm{M}_s^{\bm{h}}; \bm{C}_s^{\bm{h}})\,\dif M_s^{(i)}  + \sum_{i=1}^d \int_0^t(\mbox{$\frac{1}{2}\frac{\partial^2}{\partial x_i^2} + \frac{\partial}{\partial \Sigma_{i,i}}$})\dHe_{\bm{\alpha}}(\bm{M}_s^{\bm{h}};\bm{C}_s^{\bm{h}}) \,\dif C_s^{(i,i)} \nonumber \\
		&\quad +  \sum_{1\leq i < j\leq d} \int_0^t(\mbox{$\frac{\partial^2}{\partial x_i \partial x_j} + \frac{\partial}{\partial  \Sigma_{i,j}}$})\dHe_{\bm{\alpha}}(\bm{M}_s^{\bm{h}};\bm{C}_s^{\bm{h}}) \,\dif C_s^{(i,j)}\nonumber\\
			&=\sum_{i=1}^d \alpha_i\int_0^t \dHe_{\bm{\alpha} - \bm{e}_i}(\bm{M}_s^{\bm{h}};\bm{C}_s^{\bm{h}})\,\dif M_s^{(i)}  =\sum_{i=1}^d \alpha_i\int_0^t L_s^{\bm{\alpha} - \bm{e}_i,\bm{h}}\,\dif M_s^{(i)},
\end{align}
where the second and third sums in the first equality are equal to zero due to~\eqref{covar_pde}. Hence~$L^{\bm{\alpha},\bm{h}} $ is a continuous local martingale since it is an It\^o integral with respect to Brownian motion, see e.g.\ Proposition 2.10 in~\cite[Ch.~IV, Sec.~2]{RevuzYor}. Denoting by~$\sigma_k$ the indices of~$\bm{\alpha}$ arising from the $k$-th application of recursive relation \eqref{recursion_mvt_HMG}, one obtains the following after using~\eqref{recursion_mvt_HMG} a total of $\lvert \bm{\alpha}\rvert$ times, noting that $L_t^{(0,\ldots,0), \bm{h}} =1$: 
\begin{align*}
	L_t^{\bm{\alpha},\bm{h}} 
	&= 
	\sum_{\sigma_{\lvert \bm{\alpha}\rvert}=1}^d \alpha_{\sigma_{\lvert \bm{\alpha}\rvert}} \int_0^t L_{t_{\lvert \bm{\alpha}\rvert}}^{\bm{\alpha} - \bm{e}_{\sigma_{\lvert \bm{\alpha}\rvert}}}\,\dif M_{t_{\lvert \bm{\alpha}\rvert}}^{(\sigma_{\lvert \bm{\alpha}\rvert})} \\
&= 
\sum_{\sigma_{\lvert \bm{\alpha}\rvert}=1}^d \sum_{\sigma_{{\lvert \bm{\alpha}\rvert} -1 }=1}^d \alpha_{\sigma_{\lvert \bm{\alpha}\rvert}}(\bm{\alpha} - \bm{e}_{\sigma_{\lvert \bm{\alpha}\rvert}})_{\sigma_{{\lvert \bm{\alpha}\rvert}-1}}^+ \int_0^t \int_0^{t_{\lvert \bm{\alpha}\rvert}}  L_{t_{{\lvert \bm{\alpha}\rvert - 1}}}^{\bm{\alpha} - \bm{e}_{\sigma_{\lvert \bm{\alpha}\rvert}}- \bm{e}_{\sigma_{{\lvert \bm{\alpha}\rvert} -1} }}\,\dif M_{t_{{\lvert \bm{\alpha}\rvert}-1}}^{(\sigma_{{\lvert \bm{\alpha}\rvert} - 1})}\,\dif M_{t_{\lvert \bm{\alpha}\rvert}}^{(\sigma_{\lvert \bm{\alpha}\rvert})} \\	
	&= 
	\sum_{\sigma_{\lvert \bm{\alpha}\rvert} = 1}^d\sum_{\sigma_{{\lvert \bm{\alpha}\rvert} - 1} = 1}^d \cdots \sum_{\sigma_{1} = 1}^d \alpha_{\sigma_{{\lvert \bm{\alpha}\rvert}}} (\bm{\alpha} - \bm{e}_{\sigma_{\lvert \bm{\alpha}\rvert} } )_{\sigma_{\lvert \bm{\alpha}\rvert -1}}^+ \cdots (\bm{\alpha} - \bm{e}_{\sigma_{\lvert \bm{\alpha}\rvert}} - \cdots - \bm{e}_{\sigma_{2}})_{\sigma_{1}}^+\\
	&\quad \times \int_0^t \int_0^{t_{{\lvert \bm{\alpha}\rvert}}}\cdots \int_0^{t_2}\,\dif M_{t_1}^{(\sigma_{1})}\cdots \,\dif M_{t_{\lvert \bm{\alpha}\rvert -1}}^{(\sigma_{\lvert \bm{\alpha}\rvert - 1})}\,\dif M_{t_{\lvert \bm{\alpha}\rvert}}^{(\sigma_{\lvert \bm{\alpha}\rvert})}.
\end{align*}
It can be seen the coefficients in the sum are non-zero when the enumerations $\bm{\sigma}$ are from~$\mathcal{P}(\bm{\alpha}),$ and in that case one has 
\begin{equation*}
	\bm{\alpha}! = \alpha_{\sigma_{{\lvert \bm{\alpha}\rvert}}} (\bm{\alpha} - \bm{e}_{\sigma_{\lvert \bm{\alpha}\rvert} } )_{\sigma_{\lvert \bm{\alpha}\rvert -1}} \cdots (\bm{\alpha} - \bm{e}_{\sigma_{\lvert \bm{\alpha}\rvert}} - \cdots - \bm{e}_{\sigma_{2}})_{\sigma_{1}},
\end{equation*}
and hence representation~\eqref{multiple_integral_HMG} holds.

Using recursion relation \eqref{recursion_mvt_HMG} and It\^o's isometry one has, for $\bm{\alpha} \in \mathbb{N}^d,\ \bm{\beta} \in \mathbb{N}^d,$
\begin{equation*}
	\Ex\, L_t^{\bm{\alpha},\bm{h}}L_t^{\bm{\beta},\bm{h}} = \sum_{i=1}^d\sum_{j=1}^d\alpha_i\beta_j\int_0^t \Ex \,L_s^{\bm{\alpha} - \bm{e}_i,\bm{h}}L_s^{\bm{\beta} - \bm{e}_j,\bm{h}} \dot{h}_{i}(s)\dot{h}_{j}(s)\,\dif s.
\end{equation*} 
Repeatedly applying this recursion formula to $\Ex(L_t^{\bm{\alpha},\bm{h}})^2,$ we obtain
\begin{equation*}
	\Ex(L_t^{\bm{\alpha},\bm{h}})^2 = (\bm{\alpha}!)^2 \sum_{\bm{\sigma}, \bm{\nu} \in \mathcal{P}(\bm{\alpha})} \int_{\Delta_{\lvert \bm{\alpha} \rvert}(t)} \prod_{k = 1}^{\lvert \bm{\alpha}\rvert}  \dot{h}_{\sigma_{k}}(s_k)\dot{h}_{\nu_{k}}(s_k)\,\dif \bm{s}, \quad t\in[0,T].
\end{equation*}
Interchanging the order of summation and integral, justified by the absolute integrability of the integrand, we have
\begin{align*}
	\Ex(L_t^{\bm{\alpha},\bm{h}})^2  &=(\bm{\alpha}!)^2\int_{\Delta_{\lvert \bm{\alpha} \rvert}(t)}\sum_{\sigma,\nu \in \mathcal{P}(\bm{\alpha})}  \prod_{k = 1}^{\lvert \bm{\alpha}\rvert}   \dot{h}_{\sigma_{k}}(s_k)\dot{h}_{\nu_{k}}(s_k)\,\dif \bm{s} \nonumber \\
&=\frac{(\bm{\alpha}!)^2}{\lvert \bm{\alpha}\rvert!}\int_{[0,t]^{\lvert \bm{\alpha}\rvert}}\sum_{\sigma,\nu \in \mathcal{P}(\bm{\alpha})}  \prod_{k = 1}^{\lvert \bm{\alpha}\rvert}   \dot{h}_{\sigma_{k}}(s_k)\dot{h}_{\nu_{k}}(s_k) \,\dif \bm{s} \nonumber \\
&=\frac{(\bm{\alpha}!)^2}{\lvert \bm{\alpha}\rvert!}\sum_{\bm{\sigma}, \bm{\nu} \in \mathcal{P}(\bm{\alpha})} \prod_{k = 1}^{\lvert \bm{\alpha}\rvert}  \int_0^t  \dot{h}_{\sigma_{k}}(s) \dot{h}_{\nu_{k}}(s)\,\dif s ,
\end{align*}
where in the first equality we used the symmetry of the integrand  under permutations to extend the integration region from $\Delta_{\lvert \bm{\alpha}\rvert}(t)$ to $[0,t]^{\lvert \bm{\alpha}\rvert}.$ Using the Cauchy--Bunyakovsky--Schwarz inequality, one has
\begin{align*}\label{square_integrable_HMG}
	\Ex(L_t^{\bm{\alpha},\bm{h}} )^2 &\leq \frac{(\bm{\alpha}!)^2}{\lvert \bm{\alpha}\rvert!}\sum_{\bm{\sigma},\bm{\nu} \in \mathcal{P}(\bm{\alpha})}\prod_{k = 1}^{\lvert \bm{\alpha}\rvert} \| h_{\sigma_k}\|_{\cmSpace_t} \| h_{\nu_k}\|_{\cmSpace_t} \nonumber \\
	&=\frac{(\bm{\alpha}!)^2}{ \lvert \bm{\alpha}\rvert!  } \binom{\lvert \bm{\alpha}\rvert}{\bm{\alpha}}^2 \prod_{k=1}^{d}\| h_k\|_{\cmSpace_t}^{2\alpha_k} \nonumber \\
		&=   \lvert \bm{\alpha} \rvert! \prod_{k=1}^{d}\| h_k\|_{\cmSpace_t}^{2 \alpha_k} ,		\quad t\in [0,T],
\end{align*}
where in the second last equality we used the fact that~$\bm{\sigma}$ and~$\bm{\nu}$ are both permutations of $\bm{\alpha}$ and thus $\prod_{k=1}^{\lvert \bm{\alpha}\rvert} \| h_{\sigma_k}\|_{\cmSpace_t} = \prod_{k=1}^{\lvert \bm{\alpha}\rvert} \| h_{\nu_k}\|_{\cmSpace_t} , $ and that $\lvert \mathcal{P}(\bm{\alpha})\rvert =\binom{\lvert \bm{\alpha \rvert}}{\bm{\alpha}} = \lvert \bm{\alpha}\rvert! / \bm{\alpha}!.$ Since $\| h\|_{\cmSpace_t}  \leq \| h\|_{\cmSpace_T}$ for all $t \in [0,T],$ it follows that $\sup_{t\in[0,T]}\Ex(L_t^{\bm{\alpha},\bm{h}})^2 \leq \Ex(L_T^{\bm{\alpha},\bm{h}})^2.$ Inequality~\eqref{square_integrable_HMG} is proved.
\end{proof}

Using the auxiliary martingales~$L_t^{\bm{\alpha},\bm{h}}$ introduced in Lemma~\ref{thm:mvt_HMG}, we can now prove Theorems~\ref{thm:multi_linear_DF}--\ref{thm:psi_operator_rep}.

\begin{proof}[Proof of Theorem~\ref{thm:multi_linear_DF}]
The proof has three parts. In the first part we establish the following power series representation: for $d\geq 1,$ $\bm{z} \in \mathbb{R}^d,$ $g \in \contSpaceZero,$ $\bm{h} \in \cmSpace^d,$
\begin{equation*}
	\overline{F}[g + \bm{z}^{\top}\bm{h}](t) = \sum_{\bm{\alpha} \in \mathbb{N}_{\geq 0}^d} \frac{\bm{z}^{\bm{\alpha}} }{\bm{\alpha}!} \Ex\, L_t^{\bm{\alpha},\bm{h}}\mathbf{1}_{\tau_g > t}, \quad t \in [0,T],
\end{equation*}
and prove that $\Ex\, L_t^{\bm{1}_n,\bm{h}}\mathbf{1}_{\tau_g > t}$ is equal to the $n$-th weak $\cmSpace$-derivative of~$\overline{F}$ at~$g,$ i.e.
\begin{equation*}
\nabla_{\bm{h}}^n \overline{F}[g](t)=	\Ex\, L_t^{\bm{1}_n,\bm{h}}\mathbf{1}_{\tau_g > t} ,\quad (t,g,\bm{h}) \in [0,T] \times \contSpaceZero \times \cmSpace^n.
\end{equation*}
In the second part we prove that these weak $\cmSpace$-derivatives are in fact strong $\cmSpace$-derivatives. Lastly, we specialise the previous results to a single direction $(d =1)$ to obtain the Taylor expansion~\eqref{eqn:taylor_series} together with the remainder bound~\eqref{eqn:error_remainder_taylor}. 

We now begin the first part of the proof. By the spatial homogeneity of the Brownian motion, we  transfer the boundary perturbation~$\bm{z}^{\top}\bm{h}$ into the drift of the process to obtain
\begin{equation}\label{eqn:boundary_shift_pre_cmthm}
	\overline{F}[g + \bm{z}^{\top}\bm{h}](t) \equiv \Ex \mathbf{1}_{\tau_{g + \bm{z}^{\top}\bm{h}} > t} = \Px(W_s -  (\bm{z}^{\top}\bm{h})(s) < g(s),\ s\in [0,t]),
\end{equation}
for $t\in[0,T]$. Recalling the definition for $\bm{M}_t^{\bm{h}}$ from~\eqref{defn:vector_mg}, we now proceed to verify Novikov's criterion 
\begin{equation}\label{eqn:novikov_criterion}
\Ex\exp\{\mbox{$\frac{1}{2}$}[ \bm{z}^{\top}\bm{M}^{\bm{h}}]_T\}   <\infty,
\end{equation}
in order to apply the Cameron--Martin--Maruyama--Girsanov theorem (see e.g.\ Proposition 1.15 on p.~332 in~\cite{RevuzYor}). By the Cauchy--Bunyakovsky--Schwarz inequality, for all~$\bm{z}$ in any compact set $K \subset \mathbb{R}^d,$ one has 
\begin{align*}\label{quadratic_var_bound}
[\bm{z}^{\top}\bm{M}^{\bm{h}}]_T &=\sum_{i=1}^d\sum_{j=1}^d z_iz_j [ M^{(i)},M^{(j)} ]_T \leq \Bigg(\sum_{i=1}^d \lvert z_i \rvert  \|h_i\|_{\cmSpace} \Bigg)^2  \nonumber\\
 &\leq \Bigg(\sum_{i=1}^d \lvert z_i\rvert^2\Bigg)\Bigg( \sum_{j=1}^d \| h_j\|_{\cmSpace}^2 \Bigg) < \infty,
\end{align*}
where in the first inequality we used 
\begin{equation*}
	\lvert [ M^{(i)}, M^{(j)}]_T \rvert = \lvert \langle h_i,h_j\rangle_{\cmSpace} \rvert \leq  \|h_i\|_{\cmSpace}  \|h_j\|_{\cmSpace},
\end{equation*}
which follows from the Cauchy--Bunyakovsky--Schwarz inequality. Thus Novikov's criterion~\eqref{eqn:novikov_criterion} is satisfied and it follows that the Dol\'eans exponential 
\begin{equation*}
	\mathcal{E}(\boldsymbol{z}^{\top}\boldsymbol{M}^{\boldsymbol{h}})_t := \exp\{\boldsymbol{z}^{\top} \boldsymbol{M}^{\boldsymbol{h}}_t - \mbox{$\frac{1}{2}$}[\bm{z}^{\top}\bm{M}^{\bm{h}}]_t\}, \quad t\in [0,T],
\end{equation*}
is a uniformly integrable martingale (see e.g.\ Proposition 1.15 on p.~332 in~\cite{RevuzYor}) and hence by the Cameron--Martin--Maruyama--Girsanov theorem (see e.g.\ Proposition 1.13 on p.~331 in~\cite{RevuzYor}), one has from~\eqref{eqn:boundary_shift_pre_cmthm}
\begin{equation}\label{eqn:bcp_CMG}
	\overline{F}[g + \bm{z}^{\top}\bm{h}](t) =\Ex\,\mathcal{E}(\boldsymbol{z}^{\top}\boldsymbol{M}^{\boldsymbol{h}})_t \mathbf{1}_{\tau_g > t} , \quad t\in [0,T].
\end{equation}

Using It\^o's formula one obtains $[\bm{z}^{\top}\bm{M}^{\bm{h}} ]_t = \boldsymbol{z}^{\top}\boldsymbol{C}_t^{\bm{h}}\boldsymbol{z},$ and hence
\begin{equation*}
	\mathcal{E}(\boldsymbol{z}^{\top}\boldsymbol{M}^{\boldsymbol{h}})_t  = \exp\{\boldsymbol{z}^{\top} \boldsymbol{M}^{\boldsymbol{h}}_t - \mbox{$\frac{1}{2}$}\boldsymbol{z}^{\top}\boldsymbol{C}_t^{\bm{h}}\bm{z}\} \equiv \Lambda(\bm{z},\bm{M}_t^{\bm{h}};\bm{C}_t^{\bm{h}}), \quad t\in[0,T].
\end{equation*}
Substituting the multivariate Hermite polynomial series representation~\eqref{gen_func} for~$\Lambda$ into~\eqref{eqn:bcp_CMG}, we obtain a series representation for $\overline{F}[g+\bm{z}^{\top}\bm{h}]$ in terms of the auxiliary martingale~$L_t^{\bm{\alpha},\bm{h}} = \dHe_{\bm{\alpha}}(\bm{M}_t^{\bm{h}}; \bm{C}_t^{\bm{h}})$ introduced in~\eqref{defn:L_a_h}:
\begin{align}\label{multilinear_rep}
	\overline{F}[g + \bm{z}^{\top}\bm{h}](t) &  =  \Ex\,\Lambda(\boldsymbol{z},\boldsymbol{M}^{\boldsymbol{h}}_t ; \bm{C}_t^{\bm{h}})\mathbf{1}_{\tau_g > t}  =  \Ex  \sum_{\bm{\alpha} \in \mathbb{N}_{\geq 0}^d} \frac{\boldsymbol{z}^{\bm{\alpha}}}{\bm{\alpha}!} \dHe_{\bm{\alpha}}(\bm{M}_t^{\bm{h}} ; \bm{C}_t^{\bm{h}}) \mathbf{1}_{\tau_g > t} \nonumber\\
	&= \Ex  \sum_{\bm{\alpha} \in \mathbb{N}_{\geq 0}^d} \frac{\boldsymbol{z}^{\bm{\alpha}}}{\bm{\alpha}!} L_t^{\bm{\alpha},\bm{h}} \mathbf{1}_{\tau_g > t}, \quad t\in[0,T].
\end{align}

Next, we prove that we can interchange the order of expectation and summation in~\eqref{multilinear_rep} to obtain a uniformly convergent power series representation for~$\overline{F}$ and relate the coefficients of the power series with the weak $\cmSpace$-derivatives of~$\overline{F}.$ This interchange of order will be justified by verifying that 
\begin{equation}\label{ineq:fubini_cond}
	\sup_{t\in[0,T]}\sup_{g \in \contSpaceZero}\sum_{\bm{\alpha} \in \mathbb{N}_{\geq 0}^d} \Ex \Big \lvert\frac{\bm{z}^{\bm{\alpha}} }{\bm{\alpha}!}  L_t^{\bm{\alpha},\bm{h}} \mathbf{1}_{\tau_g > t}  \Big\rvert < \infty,
\end{equation}
and using the Fubini theorem. 

By the Cauchy--Bunyakovsky--Schwarz inequality and using the fact that $\|\overline{F}[g]\|_{\infty}$ is bounded by one, we have from inequality~\eqref{square_integrable_HMG}
\begin{align}\label{ineq:abs_ind_L_bound}
	\sup_{t\in[0,T]}\sup_{g \in \contSpaceZero}\Ex\lvert L_t^{\bm{\alpha},\bm{h}} \rvert \mathbf{1}_{\tau_g > t} &\leq \sup_{t\in[0,T]}\sup_{g \in \contSpaceZero}\big(\overline{F}[g](t)  \Ex(L_t^{\bm{\alpha,\bm{h}}} )^{2}\big)^{1/2}\leq \big(\Ex(L_T^{\bm{\alpha,\bm{h}}} )^{2}\big)^{1/2} \nonumber\\
	&\leq \sqrt{\lvert\bm{\alpha}\rvert!} \max_{1 \leq k\leq d} \| h_k \|_{\cmSpace}^{\lvert \bm{\alpha}\rvert}.
\end{align} 
Setting $c := \max_{1\leq j \leq d}\lvert z_j\rvert \max_{1\leq k \leq d}\|h_k\|_{\cmSpace},$ it follows from~\eqref{ineq:abs_ind_L_bound} that
\begin{align*}\label{abs_inf_sum}
	 \sup_{t\in[0,T]}\sup_{g \in \contSpaceZero}\sum_{\lvert \bm{\alpha}\rvert \geq 0}  \Ex\left| \frac{\boldsymbol{z}^{\bm{\alpha}}}{\bm{\alpha}!}  L_t^{\bm{\alpha},\bm{h}} \mathbf{1}_{ \tau_g > t}\right|
	&\leq  \sum_{ \lvert \bm{\alpha}\rvert \geq 0}\frac{c^{\lvert \bm{\alpha} \rvert} }{ \bm{\alpha} !}  \sqrt{\lvert\bm{\alpha}\rvert !}  
	= \sum_{n=0}^{\infty}\Bigg(\sqrt{n!}c^n\sum_{\lvert \bm{\alpha}\rvert = n} \frac{1}{ \bm{\alpha}! }\Bigg)\nonumber\\
	&=  \sum_{n = 0}^{\infty}\frac{(dc)^n}{\sqrt{n!}} < \infty,
\end{align*}
where we used the multinomial theorem  $\sum_{\lvert \bm{\alpha}\rvert =n} n!/\bm{\alpha}! =  (1+\cdots+ 1)^n =d^n $ in the second equality. Therefore~\eqref{ineq:fubini_cond} holds and thus we can use Fubini's theorem to justify the switch of order of summation and expectation in \eqref{multilinear_rep} to obtain
\begin{equation}\label{eqn:power_series}
	\overline{F}[g+ \bm{z}^{\top} \bm{h}](t) =  \sum_{\bm{\alpha} \in \mathbb{N}_{\geq 0}^d} \frac{\bm{z}^{\bm{\alpha}}}{\bm{\alpha}!} \Ex\,  L_t^{\bm{\alpha},\bm{h}} \mathbf{1}_{\tau_g > t} , \quad  t\in[0,T],
\end{equation}
for $\bm{z} \in \mathbb{R}^d,\ \bm{h} \in \cmSpace^d,$ which, for any compact subset $ K\subset \mathbb{R}^d,$ converges absolutely and uniformly for~$\bm{z} \in K.$ By the uniform convergence of the series, we can take partial derivatives with respect to each variable $z_i,$ $i =1,\ldots,d,$ of the power series representation for~$\overline{F}$ in~\eqref{eqn:power_series} and switch the order of derivative and summation, and then set~$z_i,$ $i=1,\ldots,d$ to zero, to obtain 
\begin{align}\label{eqn:derivative_is_hermite_mg}
	\nabla_{\bm{h}}^{\bm{\alpha}}\overline{F}[g](t) &:= \frac{\partial^{\lvert \bm{\alpha}\rvert}}{\partial z_1^{\alpha_1} \cdots \partial z_d^{\alpha_d}}\overline{F}\biggl[g + \sum_{i=1}^d z_i h_i\biggr]\biggr|_{z_1,\ldots,z_d = 0} \nonumber\\
	&= \Ex\, L_t^{\bm{\alpha},\bm{h}} \mathbf{1}_{\tau_g > t} , \quad (t,\bm{\alpha},\bm{h}) \in [0,T]\times  \mathbb{N}_{\geq 0}^d \times \cmSpace^d.
\end{align}

Lastly, we show that the $n$-th weak $\cmSpace$-derivative of~$\overline{F}$ is given by $\nabla_{\bm{h}}^{\bm{1}_n} \overline{F}[g]$ by verifying the $n$-linearity and boundedness (in the sense of~\eqref{defn:boundedness}) of the operator 
\begin{equation*}\label{defn:gamma}
\gamma_{n}^g[\bm{h}](t) := \nabla_{\bm{h}}^{\bm{1}_n} \overline{F}[g]=\Ex\,  L_t^{\bm{1}_n,\bm{h}} \mathbf{1}_{\tau_g > t} ,\quad (t,\bm{h}) \in [0,T] \times \cmSpace^n,
\end{equation*}
and prove that $\gamma_{n}^{g}[\bm{h}]$ is an element of $\contSpace.$ It follows that, for $\bm{h} = (h_1,\ldots,h_n) \in \cmSpace^{n},$ one has from~\eqref{multiple_integral_HMG}, 
\begin{equation*}
	L_t^{\bm{1}_n,\bm{h}} =  (-1)^n\sum_{\bm{\sigma} \in \mathcal{P}(\bm{1}_n)} \int_0^t \int_0^{t_{n}}\cdots \int_0^{t_2} \prod_{k=1}^n \dot{h}_{\sigma_k }(t_k)  \,\dif W_{t_1}\cdots\,\dif W_{t_{n-1}}\,\dif W_{t_{n}},
\end{equation*}
which is clearly $n$-linear in $\bm{h},$ and hence the $n$-linearity of~$\gamma_{n}^g$ is established. From the Cauchy--Bunyakovsky--Schwarz inequality and~\eqref{square_integrable_HMG}, since $\bm{1}_n! = 1$ and $\lvert \bm{1}_n\rvert =n,$ one obtains
\begin{equation}\label{eqn:bounded_derivative_order_n}
	\| \gamma_{n}^g[\bm{h}] \|_{\infty} \leq \sup_{t \in [0,T]}\big(\Ex (L_t^{\bm{1}_n,\bm{h}} )^2 \big)^{1/2}\leq \sqrt{n!} \prod_{k=1}^n \| h_k \|_{\cmSpace}, 
\end{equation}
establishing its boundedness in the sense of~\eqref{defn:boundedness}. We now prove that $\gamma_{n}^g[\bm{h}](t)$ is continuous in~$t \in [0,T].$ Using the fact that $L^{\mathbf{1}_n, \bm{h}}$ has continuous sample paths and~$\tau_g$ has no atoms, i.e.\ $\Px(\tau_g = t) =0$ for all $t \in [0,T],$ since~$g$ is continuous and $g(0)> 0$, we see that $\gamma_{n}^g[\bm{h}](t)$ is a continuous function in $t$:
\begin{align*}
\gamma_{n}^g[\bm{h}](t)-\gamma_{n}^g[\bm{h}](s)	&=
\lvert \Ex \, L_t^{\bm{1}_n,\bm{h}}\mathbf{1}_{\tau_g > t} -	\Ex \,L_s^{\bm{1}_n,\bm{h}}\mathbf{1}_{\tau_g > s}\rvert \\
	&\leq	 \Ex \lvert L_t^{\bm{1}_n,\bm{h}} - L_s^{\bm{1}_n,\bm{h}}\rvert +	\Ex \lvert L_s^{\bm{1}_n,\bm{h}}\rvert\lvert \mathbf{1}_{\tau_g > t} -\mathbf{1}_{\tau_g > s} \rvert\\
	&\leq \Ex \lvert L_t^{\bm{1}_n,\bm{h}} -L_s^{\bm{1}_n,\bm{h}}\rvert +	\bigl[\Ex ( L_s^{\bm{1}_n,\bm{h}})^2  \Px (\tau_g \in [s\wedge t,s\vee t))\bigr]^{1/2}\\
	&\to 0,
\end{align*}
as $s\to t$ for $t\in[0,T]$. Hence, $\gamma_{n}^g$ is bounded and $n$-linear, and $\gamma_n^g[\bm{h}]$ is an element of~$\contSpace$ for every $\bm{h} \in \cmSpace^n,$ and thus is the $n$-th weak $\cmSpace$-derivative of~$\overline{F}.$  
 
We now proceed to the second part of the proof and prove that the weak $\cmSpace$-derivatives~$\gamma_{n}^g$ are in fact strong $\cmSpace$-derivatives by showing that they are continuous with respect to~$g$ in the uniform norm and using Lemma~\ref{lemma:weak_to_strong_derivative}. Let $g_{\delta} \in \contSpaceZero  ,$ $\delta \in [0,g(0)),$ be such that $\| g_{\delta}  - g\|_{\infty}  \leq \delta.$ Recalling notation~\eqref{defn:B_g} for the sets~$B_g(t),$ it follows from the equivalence of the events
\begin{equation*}
	\{\mathbf{1}_{\tau_g > t}\} \equiv \{W \in B_g(t)\}, \quad  t \in [0,T],
\end{equation*}
and the Cauchy--Bunyakovsky--Schwarz inequality that
\begin{align*}
	&\sup_{\|h_1\|_\cmSpace=\cdots = \|h_n\|_\cmSpace = 1}\| \gamma_{n}^{g}[\bm{h}]  -\gamma_{n}^{g_{\delta}}[\bm{h}] \|_{\infty} \\
	& =\sup_{\|h_1\|_\cmSpace=\cdots= \|h_n\|_\cmSpace = 1}\sup_{t\in[0,T]}\bigl\lvert  \Ex\,L_t^{\bm{1}_n,\bm{h}} (\bm{1}_{\tau_g > t} - \bm{1}_{\tau_{g_{\delta}} > t})\bigr\rvert \\
	&\leq \bigg(\sup_{\|h_1\|_\cmSpace=\cdots =\|h_n\|_\cmSpace = 1}\Ex(L_T^{\bm{1}_n,\bm{h}})^2 \sup_{t \in[0,T]}\Px(W \in B_g(t)\Delta B_{g_{\delta}}(t)  )\bigg)^{1/2}\\
	&\leq \sqrt{n!} \sup_{t \in[0,T]}\Px(W \in B_g(t)\Delta B_{g_{\delta}}(t)  )^{1/2},
\end{align*}
where the last inequality follows from~\eqref{square_integrable_HMG}. Following the proof of Lemma~\ref{lemma:F_cont}, one sees that the right-hand side converges to zero as $\delta \to 0$ and thus $g \mapsto \gamma_{n}^g[\bm{h}] = \nabla_{\bm{h}}^{n} \overline{F}[g] $ is continuous in the operator norm, for all $n \in \mathbb{N}.$ Applying Lemma~\ref{lemma:weak_to_strong_derivative} inductively to the mappings $
g\mapsto \gamma_{n-1}^g$ for $n \in \mathbb{N}_{>0}$
 shows that each weak $\cmSpace$-derivative~$\gamma_{n}^g $ is in fact a strong $\cmSpace$-derivative of~$\overline{F}$.

Finally, we show that the Taylor expansion for $\overline{F}$ in~\eqref{eqn:taylor_series} and truncation error bound~\eqref{eqn:error_remainder_taylor} hold. By setting $\bm{z} = 1$ and $\bm{h} = h$  in~\eqref{eqn:power_series} we have
\begin{align*}
	\overline{F}[g + h](t) &= \sum_{n = 0}^{\infty} \frac{1}{n!}\Ex\,L_t^{n,h}\mathbf{1}_{\tau_g > t}= \sum_{n = 0}^{\infty} \frac{1}{n!}\Ex\,L_t^{\bm{1}_n,(h,\ldots,h)}\mathbf{1}_{\tau_g > t} \\
	&= \sum_{n = 0}^{\infty} \frac{1}{n!}\nabla_{(h,\ldots,h)}^{\mathbf{1}_n} \overline{F}[g](t) = \sum_{n = 0}^{\infty} \frac{1}{n!}\nabla_{(h,\ldots,h)}^{n} \overline{F}[g](t),
\end{align*}	
where in the second equality we used the fact that $L_t^{n,h}=L_t^{\mathbf{1}_n,(h,\ldots,h)}$ (which can be seen from the multiple integral representation~\eqref{multiple_integral_HMG}), and in the third equality we used~\eqref{eqn:derivative_is_hermite_mg}. Thus the series representation~\eqref{eqn:taylor_series} holds. From~\eqref{eqn:bounded_derivative_order_n}, it follows that the truncation error $R_m[g,h] = \sum_{n=m+1}^{\infty}\frac{1}{n!}\nabla_{(h,\ldots,h)}^n\overline{F}[g]$ has the following upper bound: 
\begin{equation*}
	 \lvert R_m[g,h](t)\rvert \leq  \sum_{n=m+1}^{\infty} \frac{1}{n!}(\Ex(L_t^{\bm{1}_n ,(h,\ldots,h) })^2 \Ex \mathbf{1}_{\tau_g > t})^{1/2} \leq \overline{F}[g](t)^{1/2}\sum_{n=m+1}^{\infty}\frac{\| h \|_{\cmSpace_t}^n}{\sqrt{n!}} 
\end{equation*}
for $t \in [0,T],$ and thus~\eqref{eqn:error_remainder_taylor} holds.
\end{proof}

\begin{proof}[Proof of Theorem~\ref{thm:chaos}]
Since~$\mathbf{1}_{\tau_g > t} \in \LpSpace^2(\Omega, \mathbf{F}^W,\Px),$ where $\mathbf{F}^W := 
\{\mathcal{F}_t^W : t \geq0\}$ denotes the natural filtration generated by~$W,$ the Wiener--It\^o chaos expansion theorem (see e.g.~\cite[Theorem 1.1.2]{Nualart2006}) states that there exists a sequence of deterministic functions~$\psi_{t}^{(n)}\in \LpSpace^2(\Delta_n(t)),$ $n\in \mathbb{N},$ uniquely determined by $\mathbf{1}_{\tau_g > t}$, such that the following mean-square convergent series representation holds: 
\begin{equation}\label{eqn:chaos_expansion_indicator}
	\mathbf{1}_{\tau_g > t}  =  \Ex\mathbf{1}_{\tau_g > t} +\sum_{n=1}^{\infty} I_{t}^{(n)}(\psi_{t}^{(n)}) ,
\end{equation}
i.e.\ \eqref{eqn:series_indicator} holds. 

We now prove~\eqref{eqn:psi_n_fla}. From~\eqref{eqn:derivative_is_hermite_mg} we have, for $\bm{h} = (h_1,\ldots,h_n) \in \cmSpace^n,$ $n \in \mathbb{N},$
\begin{equation}\label{eqn:dF_hermite_mg_rep}
	\nabla_{\bm{h}}^{\bm{1}_n} \overline{F}[g](t) = \Ex\,  L_t^{\bm{1}_n, \bm{h}} \mathbf{1}_{\tau_g > t} , \quad  t \in [0,T] .
\end{equation}	 
Denote the symmetric tensor-product of $\dot{h}_i,$ $i=1,\ldots,n,$ by
\begin{equation*}
	(\dot{h}_1 \odot \cdots \odot \dot{h}_n)(t_1,\ldots,t_n) := \frac{1}{n!}\sum_{\bm{\sigma} \in \mathcal{P}(\mathbf{1}_n) }\prod_{i=1}^n \dot{h}_{\sigma_i}(t_i),\quad (t_1,\ldots,t_n) \in (0,t)^n,
\end{equation*}
where $\mathcal{P}(\bm{1}_n)$ is the set of unique permutations of $\{1,2,\ldots,n\}.$ Using the multiple integral representation for $L_t^{\bm{1}_n,\bm{h}},$ given by~\eqref{multiple_integral_HMG} in the case when $\bm{\alpha} = \mathbf{1}_n$ and~$M^{(i)},$ $i=1,\ldots,n,$ is given by~\eqref{defn:M^h}, it follows that, for $\bm{h}\in \cmSpace^n, \ n \in \mathbb{N},$
\begin{align}\label{eqn:hermite_mg_J}
 L_t^{\bm{1}_n,\bm{h}} &= (-1)^n  \sum_{\bm{\sigma} \in \mathcal{P}(\bm{1}_n)}\int_0^{t}\int_0^{t_n}\cdots \int_0^{t_2}  \prod_{i=1}^n \dot{h}_{\sigma_i}(t_{i})\,\dif W_{t_1}\cdots \dif W_{t_{n-1}}\,\dif W_{t_n}\nonumber \\
  & = (-1)^n  I_{t}^{(n)}(\dot{h}_1 \odot \cdots \odot \dot{h}_n) ,\quad t\in[0,T].
\end{align}
Denote the symmetric extension of~$\psi_t^{(n)} $ on the simplex $\Delta_n$ to the cube $(0,t)^n$ by
\begin{equation*}
	\tilde{\psi}_{t}^{(n)}(t_1,\ldots,t_n ) := \psi_{t}^{(n)}(t_{(1)},\ldots,t_{(n)}), 
	\end{equation*}
where $(t_{(1)},\ldots,t_{(n)})$ denotes the increasing rearrangement of $(t_1,\ldots,t_n).$ Since the set of ties $t_{(i)} = t_{(j)}$ forms a set of zero Lebesgue measure and we are only interested in integrals of the function~$\tilde{\psi}_t^{(n)}$ on $(0,t)^n,$ the definition of~$\psi_t^{(n)}$ on the set of ties can be chosen arbitrarily. Using the multiple Wiener--It\^o integral representations~\eqref{eqn:chaos_expansion_indicator} and~\eqref{eqn:hermite_mg_J} for~$\mathbf{1}_{\tau_g > t}$ and~$L_t^{\bm{1}_n,\bm{h}},$ respectively, it follows from~\eqref{eqn:dF_hermite_mg_rep} and the It\^o isometry: for $m,n \in \mathbb{N}_{>0},$ $f_1,f_2 \in \LpSpace^2(\Delta_n(t)),$
\begin{equation}\label{eqn:ito_isometry}
	\Ex\, I_t^{(m)}(f_1) I_t^{(n)}(f_2) = (n!)^2 \langle f_1,f_2\rangle_{\LpSpace^2(\Delta_n(t))}\mathbf{1}_{m=n}  = n! \langle \tilde{f}_1,\tilde{f}_2\rangle_{\LpSpace^2((0,t)^n)}\mathbf{1}_{m=n}, 
\end{equation}
that, for all $n\in\mathbb{N}_{>0},$ one has 
\begin{align}\label{eqn:derivative_chaos_dual}
	\nabla_{\bm{h}}^{\mathbf{1}_n}\overline{F}[g](t) &= \overline{F}[g](t)\Ex\, I_t^{(n)}(\dot{h}_1 \odot \cdots \odot  \dot{h}_n) \nonumber\\
	&\quad + (-1)^n \Ex\, I_{t}^{(n)}(\dot{h}_1 \odot \cdots \odot  \dot{h}_n)\sum_{m=1}^{\infty} I_{t}^{(m)}(\psi_{t}^{(m)})   \nonumber\\
	&= (-1)^n n!\langle \dot{h}_1 \odot \cdots \odot \dot{h}_n,\tilde{\psi}_{t}^{(n)}\rangle_{\LpSpace^2((0,t)^n)},
\end{align}
where we used the fact that $I_t^{(n)}(\dot{h}_1 \odot \cdots \odot  \dot{h}_n)$ is a zero-mean martingale. Recalling the definition for $\kappa_{t_i}(s)\equiv s \wedge t_i,$ $s\in [0,T],$ $ (t_1,\ldots,t_n) \in \Delta_n(t),$ one has 
\begin{equation*}
(\dot{\kappa}_{t_{1}} \odot \cdots \odot  \dot{\kappa}_{t_{n}} )(s_1,\ldots,s_n) = \frac{1}{n!}\sum_{\bm{\sigma} \in \mathcal{P}(\mathbf{1}_n)}\prod_{i=1}^n \dot{\kappa}_{t_{\sigma_i}}(s_i) = \frac{1}{n!}\sum_{\bm{\sigma}\in \mathcal{P}(\mathbf{1}_n) } \prod_{i=1}^n \mathbf{1}_{s_i \leq t_{\sigma_i}},
\end{equation*}
and thus
\begin{align*}
	&\langle  \dot{\kappa}_{t_{1}} \odot \cdots \odot \dot{\kappa}_{t_{n}} , \tilde{\psi}_t^{(n)}\rangle_{\LpSpace^2((0,t)^n)} \\
	&=
\frac{1}{n!}\sum_{\bm{\sigma} \in \mathcal{P}(\mathbf{1}_n)}	\int_0^{t_{\sigma_n}}\int_0^{t_{\sigma_{n-1}}} \cdots \int_0^{t_{\sigma_1} }  \tilde{\psi}_t^{(n)}(s_1,\ldots,s_n)\,\dif s_1\cdots \dif  s_{n-1}\,\dif s_n.
\end{align*}
Since $\tilde{\psi}_t^{(n)} \in \LpSpace^2((0,t)^n),$ and thus absolutely integrable, each integral in the sum on the right-hand side is absolutely continuous in each variable $t_{\sigma_1},\ldots,t_{\sigma_n}$ separately (for a.e.\ fixed values of the remaining variables). By iterated application of the one-dimensional Lebesgue differentiation theorem, the following mixed derivative exists for a.e.\ $(t_1,\ldots,t_n) \in \Delta_n(t)$:
\begin{align*}
	\frac{\partial^n}{\partial t_1 \cdots \partial t_n}  \langle   \dot{\kappa}_{t_{1}} \odot \cdots \odot \dot{\kappa}_{t_{n}} , \tilde{\psi}_{t}^{(n)}\rangle_{\LpSpace^2((0,t)^n)} &= \frac{1}{n!}\sum_{\bm{\sigma} \in \mathcal{P}(\mathbf{1}_n)}\tilde{\psi}_t^{(n)} (t_{\sigma_1},\ldots,t_{\sigma_n}) \\
	&=\psi_{t}^{(n)}(t_1,\ldots,t_n).
\end{align*}
Combining this with~\eqref{eqn:derivative_chaos_dual} in the case $\bm{h} = (\kappa_{t_1},\ldots,\kappa_{t_n})$ we obtain the claimed representation~\eqref{eqn:psi_n_fla}. 

Using~\eqref{eqn:derivative_chaos_dual} with $\bm{h} = (h,\ldots,h)$ and using the symmetry of $(\dot{h} \odot \cdots \odot \dot{h})(t_1,\ldots,t_n) = \prod_{i=1}^n \dot{h}(t_i)$ and $\tilde{\psi}_t^{(n)}$ in $(0,t)^n,$ one has
\begin{align*}
	\nabla_{\bm{h}}^{\mathbf{1}_n}\overline{F}[g](t) &=\nabla_{(h,\ldots,h)}^{n}\overline{F}[g](t) \\
	&= (-1)^nn!\langle \dot{h} \odot \cdots \odot \dot{h},\tilde{\psi}_{t}^{(n)}\rangle_{\LpSpace^2((0,t)^n)}\\
	&= (-1)^n (n!)^2\int_{0}^{t}\int_{0}^{t_n} \cdots\int_{0}^{t_2} \psi_{t}^{(n)}(t_1,\ldots,t_n)\prod_{i=1}^n \dot{h}(t_i)\,\dif t_1\,\cdots\,\dif t_{n-1}\,\dif t_n,
\end{align*}
which is the claimed~\eqref{eqn:dF_multiple_integral}.
\end{proof}

\begin{remark}
We note that no specific properties of the first-passage time were used in the proofs of Theorems~\ref{thm:multi_linear_DF} and~\ref{thm:chaos} except in the proof of the strong $\cmSpace$-differentiability of~$\overline{F}$. Thus analogous results to Theorems~\ref{thm:multi_linear_DF} and~\ref{thm:chaos} for general expected values of square-integrable Brownian functionals as functions of time-dependent drift functions can also be obtained using the same proof technique. 
\end{remark}

The following martingale representation is a generalisation of the one in Lemma~7 from~\cite{Liang2025} to the case when the terminal function~$\phi$ is not identically equal to one.

\begin{lemma}\label{lemma:ito_fla_stop}
Suppose that $g \in \contSpaceZero \cap\contSpace^2$ and $t \in(0,T].$ If $\phi(t,\,\cdot\,)$ is bounded and continuous on $(-\infty,g(t)),$ then one has
\begin{equation}\label{eqn:ito_fla_phi}
	\phi(t,W_t)\mathbf{1}_{\tau_g > t} = (Q_t\phi)(0,0) + \int_0^t  (\Dx Q_t \phi)(s, W_s)\mathbf{1}_{\tau_g > s}\,\dif W_s \quad  a.s.
\end{equation}
\end{lemma}

\begin{proof}
For $t\in (0,T],$ let
\begin{equation*}
	w_t(s,x) := (Q_t\phi)(s,x) \equiv \Ex\,\phi(t,W_t^{s,x})\mathbf{1}_{\tau_g^{s,x} > t}, \quad (s,x) \in D_t.
\end{equation*} 
By replacing the terminal condition~$\bm{1}_{s = T}$ in the terminal-boundary-value problem (4.3) from~\cite{Liang2025} with $\phi(t,x)\mathbf{1}_{s=t}$ and following the proof of Lemma 1 in~\cite{Liang2025}, one shows that $w_t  \in \contSpace^{1,2}(D_t) \cap \contSpace(\text{cl}(D_{t-\varepsilon}))$ for all $\varepsilon \in (0,t)$ and is a bounded solution to the following terminal-boundary-value problem for the backward Kolmogorov equation:
\begin{equation}\label{eqn:BKE}
\begin{cases}
	(\mbox{$\frac{\partial}{\partial s}$}   + \mbox{$\frac{1}{2}$}\mbox{$\frac{\partial^2}{\partial x^2}$})w_t(s,x) =0, & (s,x) \in D_t,\\
	\lim_{s\uparrow t}w_t(s,x) = \phi(t,x), & x < g(t),\\
	\lim_{x\uparrow g(s)}w_t(s,x) =0 ,& s\in[0,t).
	\end{cases} 
\end{equation}
Therefore, one has from It\^o's formula,
\begin{equation}\label{eqn:ito_fla_stopped}
w_t(s \wedge \tau_g , W_{ s \wedge \tau_g })= w_t(0,0)+\int_0^{ s \wedge \tau_g }\Dx w_t(u, W_u )\,\dif W_u, \quad s\in[0,t).
\end{equation}
Since $w_t(\tau_g ,W_{\tau_g  } ) = w_t(\tau_g,g(\tau_g)) = 0$ a.s.\ on the event $\{\tau_g  \leq  t\},$ the left-hand side of~\eqref{eqn:ito_fla_stopped} can be expressed as  
\begin{equation}\label{eqn:boundary_cond}
w_t( s \wedge \tau_g  ,W_{ s \wedge \tau_g } ) = w_t(s,W_s )\mathbf{1}_{\tau_g >  s} + w_t(\tau_g , W_{\tau_g } )\mathbf{1}_{\tau_g  \leq s}  = w_t(s,W_s )\mathbf{1}_{\tau_g >  s}.
\end{equation}
Combining~\eqref{eqn:ito_fla_stopped} and~\eqref{eqn:boundary_cond}, the following holds for all $s \in [0,t)$: 
\begin{equation}\label{eqn:operator_recursion}
	M_s:= w_t(s,W_s )\mathbf{1}_{\tau_g > s} = w_t(0,0) + \int_0^s  (\Dx w_t)(u,W_u)\mathbf{1}_{\tau_g > u}\,\dif W_u.
\end{equation}

By the terminal condition from~\eqref{eqn:BKE}, the left-hand side of~\eqref{eqn:operator_recursion} converges to $M_t= \phi(t,W_t)\mathbf{1}_{\tau_g > t}$ almost surely as $s\uparrow t$ since $\Px(\tau_g = t )=0.$ 

Since $w_t$ is bounded on $D_t,$ it follows that~$\sup_{s \in [0,t)}\lvert M_s\rvert \leq \sup_{(s,x)\in D_t}\lvert w_t(s,x)\rvert$ is also bounded, and therefore the local martingale on the right-hand side of~\eqref{eqn:operator_recursion} is a square-integrable martingale on $[0,t)$.  By applying It\^o's isometry, one obtains from~\eqref{eqn:operator_recursion}
\begin{equation*}
	\Ex M_s^2 = w_t(0,0)^2 + \Ex \int_0^s (\Dx w_t)(u,W_u)^2 \mathbf{1}_{\tau_g > u}\,\dif u.
\end{equation*}
Passing to the limit as $s \uparrow t,$ monotone convergence and $\Ex M_t^2 \leq \sup_{x < g(t)}\lvert \phi(t,x)\rvert^2 <\infty,$ yields
\begin{equation*}
	\Ex \int_0^t (\Dx w_t)(u,W_u)^2 \mathbf{1}_{\tau_g > u}\,\dif u = \Ex M_t^2 - 	w_t(0,0)^2 < \infty.
\end{equation*}
Hence the stochastic integral on the right-hand side of~\eqref{eqn:operator_recursion} converges in $\LpSpace^2$ as $s\uparrow t.$ Thus the claim~\eqref{eqn:ito_fla_phi} follows from~\eqref{eqn:operator_recursion} as $s\uparrow t$.
\end{proof}

\begin{lemma}\label{lemma:gaussian_bound}
Suppose that $g \in \contSpaceZero \cap \contSpace^2$ and let $0<s<t\leq T.$ There exists a positive constant~$C$ depending on~$t$ and~$g$ such that
\begin{equation*}
	\sup_{x < g(s)}\lvert (\Dx Q_t \phi)(s,x) \rvert \leq \frac{C}{\sqrt{t-s}}\sup_{x < g(t)}\lvert \phi(t,x)\rvert.
\end{equation*}
\end{lemma}
\begin{proof}
The claim follows from the proof of Lemma 2 in~\cite{Liang2025}. By using the Gaussian bounds on the first spatial derivative of the taboo transition density one immediately obtains the claimed result.
\end{proof}

\begin{proof}[Proof of Theorem~\ref{thm:psi_operator_rep}]
We will prove the claimed identity by testing against the class of smooth functions whose support is compactly contained in the open simplex $\Delta_n(t),$ which we denote by $\contSpace_c^{\infty}(\Delta_n(t)).$  Multiplying the chaos expansion for $\mathbf{1}_{\tau_g > t}$ in~\eqref{eqn:series_indicator} by $I_t^{(n)}(f),$ $f \in \contSpace_c^{\infty}(\Delta_n(t)),$ and using It\^o's isometry~\eqref{eqn:ito_isometry}, one has
\begin{equation}\label{eqn:test_function}
	\Ex  \mathbf{1}_{\tau_g > t }I_t^{(n)}(f)= (n!)^2 \int_{\Delta_n(t)}\psi_t^{(n)}(\bm{t})f(\bm{t})\,\dif \bm{t}.
\end{equation}
We will now obtain another integral expression for the left-hand side by using a recursive representation for $I_t^{(n)}(f),$ which will require the following notation. Let $n \in \mathbb{N}_{\geq 1}.$ For $0 \leq k \leq n-1,$ $ s < t_{k+2} <\cdots<t_n,$ define
\begin{equation*}
	J_{k}^f(s,t_{k+2},\ldots,t_n) := \int_{\Delta_k(s)}f(t_1,\ldots,t_k,s,t_{k+2},\ldots,t_n)\,\dif W_{t_1}\cdots \,\dif W_{t_k},
\end{equation*}
with $J_0^f(s,t_2,\ldots,t_n) := f(s,t_2,\ldots,t_n).$ This satisfies the recursion relation
\begin{equation}\label{eqn:J_krecrusion}
	J_{k}^f(s,t_{k+2},\ldots,t_n) = \int_0^{s} J_{k-1}^f(t_k,s,t_{k+2},\ldots, t_n) \,\dif W_{t_k}, \quad k \geq 1.
\end{equation}
  Note that by It\^o's isometry one has
\begin{equation}\label{eqn:J_moment}
		\Ex  J_{k}^f(s,t_{k+2},\ldots,t_n)^2 = \int_{\Delta_k(s)}f(t_1,\ldots,t_k,s,t_{k+2},\ldots,t_n)^2 \,\dif t_1\cdots \dif t_k.
\end{equation}

It follows from Lemma~\ref{lemma:ito_fla_stop} in the case when $\phi \equiv 1$ that, for $t \in(0,T],$
\begin{equation}\label{eqn:MG_rep_ind}
	\mathbf{1}_{\tau_g > t} = (Q_t1)(0,0) + \int_0^t (\Dx Q_t 1)(s,W_s)\mathbf{1}_{\tau_g > s}\,\dif W_s  .
\end{equation}
Multiplying both sides by $I_t^{(n)}(f)= n!\int_0^t J_{n-1}^f(s)\,\dif W_s$ and taking expected values, we get from It\^o's isometry and~\eqref{eqn:MG_rep_ind} that
\begin{align}\label{eqn:E1I_n}
	\Ex \mathbf{1}_{\tau_g > t }I_t^{(n)}(f) &= n! \Ex  \int_0^t (\Dx Q_t1)(t_n,W_{t_n})\mathbf{1}_{\tau_g > t_n} J_{n-1}^f(t_n)\,\dif t_n.
\end{align}
To interchange the order of expectation and integration on the right-hand side using Fubini's theorem, we will verify the absolute integrability of the integrand. 

Let $t_{n+1}:=t$ and $t_{n+2}=T$ and suppose that~$\phi$ satisfies the conditions of Lemma~\ref{lemma:ito_fla_stop}. Using  the Cauchy--Bunyakovsky--Schwarz inequality twice, the It\^o integral representation \eqref{eqn:ito_fla_phi}, It\^o's isometry and~\eqref{eqn:J_moment}, we obtain the following inequality for $k=1,\ldots,n$:
\begin{align}\label{ineq:cross_moment_bound}
&\int_0^{t_{k+1}} \Ex \bigl\lvert (\Dx Q_{t_{k+1}} \phi)(s,W_s)\mathbf{1}_{\tau_g > s} J_{k-1}^f(s,t_{k+2},\ldots,t_n)\bigr\rvert\,\dif s \nonumber\\
&\leq \int_0^{t_{k+1}} \bigl(\Ex (\Dx Q_{t_{k+1}} \phi)(s,W_s)^2\mathbf{1}_{\tau_g > s} \bigr)^{1/2} \bigl(\Ex J_{k-1}^f(s,t_{k+2},\ldots,t_n)^2\bigr)^{1/2}\,\dif s \nonumber\\
&\leq \biggl(\int_0^{t_{k+1}} \Ex  (\Dx Q_{t_{k+1}} \phi)(s,W_s)^2\mathbf{1}_{\tau_g > s}\,\dif s\biggr)^{1/2}
\biggl(\int_0^{t_{k+1}} \Ex  J_{k-1}^f(s,t_{k+2},\ldots,t_n)^2 \,\dif s\biggr)^{1/2} \nonumber\\
&= \var\,(\phi(t_{k+1},W_{t_{k+1}})\mathbf{1}_{\tau_g > t_{k+1}})^{1/2}
\| f(\,\cdot\,,t_{k+1},\ldots,t_n)\|_{\LpSpace^2(\Delta_{k}(t_{k+1}))}\nonumber\\
&\leq \sup_{x <g(t_{k+1})}\lvert \phi(t_{k+1},x)\rvert\,\| f(\,\cdot\,,t_{k+1},\ldots,t_n)\|_{\LpSpace^2(\Delta_{k}(t_{k+1}))},
\end{align}
which is bounded for $t_{k+1} \in (0,t_{k+2})$ since the right-hand vanishes for $t_{k+1}$ sufficiently close to $t_{k+2}$ because the support of~$f$ is compactly contained in the open simplex.

It follows from~\eqref{ineq:cross_moment_bound}, in the case when $k = n,$ that we can exchange the order of expectation and integral in~\eqref{eqn:E1I_n} using the Fubini theorem to obtain 
\begin{equation}\label{ito_isom_f}
\Ex \mathbf{1}_{\tau_g > t }I_t^{(n)}(f) =  n! \int_0^t \Ex  (\Dx Q_t1)(t_n,W_{t_n})\mathbf{1}_{\tau_g > t_n} J_{n-1}^f(t_n)\,\dif t_n .
\end{equation}

Applying Lemma~\ref{lemma:ito_fla_stop} to the case when $\phi(s,x) = (\Dx Q_t 1)(s,x)$ (where the conditions of the lemma are satisfied due to Lemma~\ref{lemma:gaussian_bound} and the smoothness of the taboo transition density), using the recursion~\eqref{eqn:J_krecrusion}, taking expected values, and applying It\^o's isometry together with the fact that $\Ex  \,J_{n-1}^f(t_n) =0$, one has
\begin{align*}
	&\Ex  (\Dx Q_t1)(t_n,W_{t_n})\mathbf{1}_{\tau_g > t_n} J_{n-1}^f(t_n)  \\
	&=  \Ex\int_0^{t_n} (\Dx Q_{t_n} \Dx Q_{t}1)(t_{n-1},W_{t_{n-1}}) \mathbf{1}_{\tau_g > t_{n-1}} J_{n-2}^f(t_{n-1},t_n)\,\dif t_{n-1}\\
		&=  \int_0^{t_n} \Ex(\Dx Q_{t_n} \Dx Q_{t}1)(t_{n-1},W_{t_{n-1}}) \mathbf{1}_{\tau_g > t_{n-1}} J_{n-2}^f(t_{n-1},t_n)\,\dif t_{n-1},
\end{align*}
where the exchange in the order of expectation and integral is justified by the Fubini theorem, since, by the Cauchy--Bunyakovsky--Schwarz inequality and~\eqref{ineq:cross_moment_bound} in the case $k=n-1$ and Lemma~\ref{lemma:gaussian_bound}, there exists a constant $C >0$ such that
\begin{align*}
	&\int_0^{t_n}\Ex  \lvert (\Dx Q_{t_n} \Dx Q_{t}1)(t_{n-1},W_{t_{n-1}}) \mathbf{1}_{\tau_g > t_{n-1}} J_{n-2}^f(t_{n-1},t_n)\rvert \,\dif t_{n-1} \\
	&\leq \sup_{x <g(t_{n})}\lvert (\Dx Q_t1)(t_{n},x)\rvert\| f(\,\cdot\,,t_n)\|_{\LpSpace^2(\Delta_{n-1}(t_{n}))} \\
	&\leq \frac{C}{\sqrt{t-t_{n}}}\| f(\,\cdot\,,t_n)\|_{\LpSpace^2(\Delta_{n-1}(t_{n}))},
\end{align*}
which is bounded for all~$t_n \in (0,t)$ due to the compact support of $f.$

For $n \geq 2,$ repeating the above argument $n-1$ times and retaining the constant term in the final step $k=0$ (where $J_0^f(t_1,\ldots,t_n) = f(t_1,\ldots,t_n)$ is deterministic and has nonzero expectation) yields
\begin{equation}\label{eqn:target_limit}
	\Ex  \mathbf{1}_{\tau_g > t } I_t^{(n)}(f ) = n! \int_{\Delta_n(t)}f(\bm{t})A_n(\bm{t})\,\dif \bm{t},
\end{equation}
where, for $(t_1,\ldots,t_n) \in \Delta_n(t),$
\begin{equation*}
	A_n(t_1,\ldots,t_n) := (Q_{t_1}\Dx Q_{t_2}\cdots \Dx Q_{t_n}\Dx Q_t 1)(0,0),
\end{equation*}
which is locally integrable on the simplex $\Delta_n(t)$ since there exists a constant $C >0$ such that 
\begin{equation*}
	\int_{\Delta_n(t)}\lvert A_n(\bm{t})\rvert\,\dif \bm{t} \leq C\int_{\Delta_n(t)}\prod_{i=1}^{n+1} (t_{i} - t_{i-1})^{-1/2}\,\dif t_1\cdots \dif t_n = Ct^{(n-1)/2}\frac{\Gamma(\frac{1}{2})^{n+1}}{\Gamma(\frac{n+1}{2})} < \infty,
\end{equation*}	
where the inequality follows from repeatedly applying Lemma~\ref{lemma:gaussian_bound}. The case $n=1$ is immediate from~\eqref{ito_isom_f}. After subtracting~\eqref{eqn:target_limit} from~\eqref{eqn:test_function} and dividing by $n!$, we have for every $f \in \contSpace_c^{\infty}(\Delta_n(t)),$
\begin{equation*}
	\int_{\Delta_n(t)} f(\bm{t})\biggl(\psi_t^{(n)}(\bm{t}) - \frac{1}{n!}A_n(\bm{t})\biggr)\,\dif \bm{t} =0.
\end{equation*}
By the du Bois-Reymond lemma (see e.g.\ Corollary 4.24 in~\cite{Brezis2011}), it follows that~$\psi_t^{(n)}$ is equal to $\frac{1}{n!}A_n$ almost everywhere in~$\Delta_n(t)$ as claimed.
\end{proof}
 
 \subsection{The first and second $\cmSpace$-derivative of~$\overline{F}$}

 For notational convenience, for $t \in (0,T],$ let
\begin{equation}\label{defn:v_t}
	v_t(s,x) := (Q_t1)(s,x) = \Px (\tau_g^{s,x} > t), \quad (s,x) \in D_t.
\end{equation} 
We note that the existence of the spatial derivatives
\begin{equation*}
	v_t'(s,x) :=  \frac{\partial}{\partial x} v_t(s,x),\quad v_t''(s,x) := \frac{\partial^2}{\partial x^2}v_t(s,x), \quad v_t'''(s,x) := \frac{\partial^3}{\partial x^3}v_t	(s,x)
\end{equation*}
is guaranteed in $(s,x) \in D_t$ under the assumption that $g\in\contSpace^2$ by Lemma 3 in~\cite{Liang2025}. The following novel identity will be crucial in the proofs of Propositions~\ref{prop:deriv_bm} and~\ref{prop:deriv_2_bm}. 
\begin{lemma}\label{lemma:heat_gradient_identity}
Suppose that $g \in \contSpaceZero\cap\contSpace^2,$ $t\in (0,T],$ and set $f_{\tau_g^{s,x}}(u) := \Px (\tau_g^{s,x} \in \dif u)/\dif u,$ $u \in (0,t),$ $(s,x)\in D_t.$ Then, for $0 \leq s < r < t \leq T,$ $(s,x)\in D_t,$ one has
\begin{equation}\label{eqn:dxQphi}
	   (Q_r v_t')(s,x) =   \int_r^t v_t'(u,g(u)-)  f_{\tau_g^{s,x}}(u)\,\dif u.
\end{equation}
\end{lemma}

\begin{proof}
To prove the claim of the lemma, we will use the Doob $h$-transform method from the proof of Proposition~1 in~\cite{Liang2025} to obtain a probabilistic expression for $(Q_rv_t')(s,x),$ from which we can obtain~\eqref{eqn:dxQphi}.

From Lemma~3 in~\cite{Liang2025} we know that $v_t' \in \contSpace^{1,2}(D_t) \cap \contSpace(\text{cl}({D}_{t-\varepsilon}))$  is bounded in $D_{t-\varepsilon}$ for $\varepsilon \in (0,t).$ Since $v_t(s,x)$ satisfies the time-reversed heat equation $(\frac{\partial}{\partial s} + \frac{1}{2}\frac{\partial^2}{\partial x^2})v_t(s,x) = 0$ for $(s,x) \in D_t$ (Lemma 1 in~\cite{Liang2025}), it follows by linearity of differentiation and the differential operator~$\frac{\partial}{\partial s} + \frac{1}{2}\frac{\partial^2}{\partial x^2}$  that~$v_t'$ satisfies $(\frac{\partial}{\partial s} + \frac{1}{2}\frac{\partial^2}{\partial x^2})v_t'(s,x) = 0$ for $(s,x) \in D_t$. Thus, by It\^o's formula, the stopped process
\begin{equation}\label{defn:N}
	N_u^{s,x} := \frac{v_t'(u \wedge \tau_g^{s,x} ,W_{u \wedge \tau_g^{s,x}}^{s,x})}{v_t'(s,x)} ,\quad (s,x)\in D_t,\ u \in (s,t),
\end{equation}
satisfies the linear stochastic equation
\begin{equation*}
	N_u^{s,x} = 1 + \int_s^{u \wedge \tau_g^{s,x}} \Big(\frac{v_t''}{v_t'}\Big)(r,W_r^{s,x})N_r^{s,x}\,\dif W_r, \quad (s,x)\in D_t,\ u\in (s,t),
\end{equation*}
where we used the fact that~$v_t'(s,x)$ is non-zero and does not change sign in $(s,x) \in D_t$ so that one can write $v_t''(r,W_r^{s,x}) = (v_t''/v_t')(r,W_r^{s,x})N_r^{s,x}$ on the event $\{r < \tau_g^{s,x}\}$. The solution to the stochastic equation is given by the stochastic exponential 
\begin{equation*}
	N_u^{s,x} = \exp\{V_u^{s,x} - \mbox{$\frac{1}{2}$}[V^{s,x}]_u\}, \quad V_u^{s,x} := \int_s^{u \wedge \tau_g^{s,x}} \Big(\frac{v_t''}{v_t'}\Big)(r,W_r^{s,x}) \,\dif W_r
\end{equation*}
for $(s,x) \in D_t,$ $u\in (s,t).$ Since $\lvert v_t'(s,x)\rvert$ is bounded for $(s,x)\in \text{cl}(D_{t-\varepsilon})$ for any $\varepsilon \in (0,t),$ it follows that $\{N_{u\wedge (t - \varepsilon)}^{s,x} : u \in [0,t]\}$ is a bounded (positive) martingale. 

For $\varepsilon \in(0,t),$ denote by $\mathbf{Q}_t^{\varepsilon}$ the measure on $\mathcal{F}_t^W$ defined by $
\dif \mathbf{Q}_t^{\varepsilon}/\dif \Px   = N_{t-\varepsilon}^{s,x},$ $(s,x) \in D_t.$ Using the same argument as in the proof of Proposition 1 in~\cite{Liang2025}, we can pass to the limit as $\varepsilon \downarrow 0$ to obtain the existence of a probability measure~$\mathbf{Q}_t$ on $\sigma(\bigcup_{s<t}\mathcal{F}_s^W)$ such that
\begin{equation}\label{eqn:change_of_measure}
	\mathbf{Q}_t(A) = \Ex\,N_t^{s,x}\mathbf{1}_{A}, \quad A\in \mathcal{F}_t^W,\ (s,x) \in D_t. 
\end{equation}
Hence, for $ (s,x) \in D_t,$ $r \in (s,t),$ we have from~\eqref{defn:N} and~\eqref{eqn:change_of_measure} that
\begin{equation}\label{eqn:Qphi_k}
(Q_{r}v_t')(s,x)\equiv  \Ex\,v_t'(r,W_r^{s,x})\mathbf{1}_{\tau_g^{s,x} > r} = v_t'(s,x) \mathbf{Q}_t(\tau_g^{s,x} > r).
\end{equation}

With the goal of deriving an expression for $\mathbf{Q}_t(\tau_g^{s,x} > r),$ we now prove that 
\begin{equation}\label{eqn:tau_support_Q}
	\mathbf{Q}_t(\tau_g^{s,x} \leq t) = 1, \quad (s,x) \in D_t.
\end{equation}
Dividing both sides of~\eqref{eqn:Qphi_k} by $v_t'(s,x),$ using the integral representation~\eqref{eqn:operator_taboo_transition} for~$Q_t$ in terms of the taboo transition density~$q,$ and applying integration by parts and passing to the limit as $r\uparrow t$, one has
\begin{align*}	\lim_{r\uparrow t}\mathbf{Q}_t(\tau_g^{s,x} > r) &= \lim_{r\uparrow t}\frac{1}{v_t'(s,x)} \int_{-\infty}^{g(r)}v_t'(r,y)q(s,x;r,y)\,\dif y \\
	&=  -\lim_{r\uparrow t}\frac{1}{v_t'(s,x)} \int_{-\infty}^{g(r)}v_t(r,y) \mbox{$\frac{\partial}{\partial y}$}q(s,x;r,y)\,\dif y\\
	&= -\frac{1}{v_t'(s,x)} \int_{-\infty}^{g(t)}  \mbox{$\frac{\partial}{\partial y}$}q(s,x;t,y)\,\dif y\\
	&=0,
\end{align*}
where we used the fact that $\lim_{r \to t}v_t(r,y) =1$ for $y < g(t)$ in the second last equality and the boundary conditions
\begin{equation}\label{eqn:bc_taboo}
\lim_{y \uparrow g(u)}  q(s,x;u,y) = \lim_{y \downarrow -\infty} q(s,x ; u,y)= 0, \quad u \in(s,t),
\end{equation} 
in the last equality. The interchange in order of limit and integral in the second last equality is justified by the dominated convergence theorem (see the proof of Theorem~1 in~\cite{Liang2025}). This proves~\eqref{eqn:tau_support_Q}. 

Denote the density of $\tau_g^{s,x}$ under $\Px$ and $\mathbf{Q}_t$ by $f_{\tau_g^{s,x}}$ and $f_{\tau_g^{s,x}}^{\mathbf{Q}_t},$ respectively, where their existence is guaranteed by Theorem~1 and Corollary~1 in~\cite{Liang2025}.  Using~\eqref{eqn:tau_support_Q} and Corollary~1 in~\cite{Liang2025} for the representation for $f_{\tau_g^{s,x}}^{\mathbf{Q}_t}(u) = \frac{v_t'(u,g(u) - )}{v_t'(s,x)}f_{\tau_g^{s,x}}(u),$ we have, for $(s,x)\in D_t,$ $r \in(s, t),$
\begin{align}\label{eqn:Q_CDF}
	\mathbf{Q}_t(\tau_g^{s,x} > r) &= \int_r^t f_{\tau_g^{s,x}}^{\mathbf{Q}_t}(u)\,\dif u = \int_r^t\frac{v_t'(u,g(u)-)}{v_t'(s,x)}f_{\tau_g^{s,x}}(u)\,\dif u .
\end{align}
Multiplying both sides by $v_t'(s,x),$ identifying the resulting left-hand side with $(Q_rv_t')(s,x)$ using~\eqref{eqn:Qphi_k}, it follows from~\eqref{eqn:Q_CDF} that~\eqref{eqn:dxQphi} holds. 

\end{proof}
 Using the results from Theorems~\ref{thm:chaos} and~\ref{thm:psi_operator_rep} we now prove Propositions~\ref{prop:deriv_bm}--\ref{prop:approx_m}. 
  
\begin{proof}[Proof of Proposition~\ref{prop:deriv_bm}]

Using~\eqref{eqn:dF_multiple_integral} from Theorem~\ref{thm:chaos} in the case $n=1,$ the first $\cmSpace$-derivative of~$\overline{F}$ can be expressed as 
\begin{equation}\label{eqn:D_h}
	\nabla_h \overline{F}[g](t) = -\int_0^t \dot{h}(s)\psi_{t}^{(1)}(s)\,\dif s ,  \quad t \in (0,T],
\end{equation}
where it follows from~\eqref{eqn:equiv_psi} from Theorem~\ref{thm:psi_operator_rep} in the case $n=1$ and~\eqref{eqn:dxQphi} from Lemma~\ref{lemma:heat_gradient_identity} that the claim~\eqref{defn:first_chaos_linear} holds, i.e.\
\begin{equation*}
	\psi_t^{(1)}(s) = (Q_s \Dx Q_t1)(0,0) = \int_s^t v_t'(u,g(u)-)f_{\tau_g}(u)\,\dif u, \quad s \in (0,t).
\end{equation*}
Applying integration by parts to~\eqref{eqn:D_h} and using the boundary conditions $h(0) = \lim_{s\to t}\psi_t^{(1)}(s)= 0,$ we get~\eqref{eqn:dF_1_BM}. 

In the special case of the linear boundary $g(t) = a + bt,$ $a >0,$ $b \in \mathbb{R},$ one obtains the formula for $f_{\tau_g}(t)$ from the Bachelier--L\'evy formula (see e.g.\ 2.0.2 in~\cite[p.~204]{Borodin2015}). In the linear boundary case one also has, for $x < g(s),$ $0 \leq s < t,$ 
\begin{equation}\label{eqn:v_absx}
v_t(s,x) =  \Phi_{t-s}(a - x+bt) -e^{-2(a+bs-x)b}\Phi_{t-s}(-a -bs +x + b(t-s)),
\end{equation}
see e.g.\ 1.1.4 on p.\ 256 in~\cite{Borodin2015}. By taking partial derivatives with respect to~$x$ and passing to the limit as $x\uparrow g(s)$ one obtains~\eqref{eqn:bdry_grad_bm}. 
\end{proof}

\begin{proof}[Proof of Proposition~\ref{prop:deriv_2_bm}]
Let
\begin{equation*}
	J_t(t_2) := \int_0^{t_2}\dot{h}(t_1)\psi_{t}^{(2)}(t_1,t_2)\,\dif t_1, \quad 0< t_2 < t\leq T.
\end{equation*}
Using~\eqref{eqn:dF_multiple_integral} in the case $n=2$ and applying integration by parts, we have
\begin{align}\label{eqn:preIBP}
\nabla_{(h,h)}^2 \overline{F}[g](t)&=4\int_0^t\int_0^{t_2} \dot{h}(t_1)\dot{h}(t_2) \psi_{t}^{(2)}(t_1,t_2)\,\dif t_1\,\dif t_2 = 4\int_0^t \dot{h}(t_2) J_t(t_2)\,\dif t_2 \nonumber\\
&= 4h(\,\cdot\,)J_t(\,\cdot\,)|_{0}^t - 4\int_0^t h(t_2)\mbox{$\frac{\partial}{\partial t_2}$}J_t(t_2)\,\dif t_2,
\end{align}
where the boundary terms are equal to zero since $h(0)=0$ because $h \in \cmSpace,$ and $\psi_{t}^{(2)}(t_1,t-)=0$ which we will see from~\eqref{eqn:psi_2} below. Using the Leibniz integral rule,  one has
\begin{equation}\label{eqn:dI}
	\mbox{$\frac{\partial}{\partial t_2}$} J_t(t_2) = \dot{h}(t_2) \psi_{t}^{(2)}(t_2-,t_2) + \int_0^{t_2} \dot{h}(t_1)\mbox{$\frac{\partial}{\partial t_2} $}\psi_{t}^{(2)}(t_1,t_2)\,\dif t_1,
\end{equation}
where we will see that~$\psi_{t}^{(2)}(t_2-,t_2)$ is finite  and integrable with respect to $t_2 \in (0,t),$ from~\eqref{eqn:psi_2} below because
\begin{equation}\label{eqn:K_ru_diagonal}
	K(r, s) = -2\frac{a}{s}\varphi_s(a+bs)\bigg( \frac{1}{\sqrt{2 \pi (s-r)}} + O(1)\bigg),\quad r \uparrow s,
\end{equation}
which we will show using its closed-form expression~\eqref{eqn:fla_K_ru_explicit} below. Substituting~\eqref{eqn:dI} into~\eqref{eqn:preIBP} proves the claimed expression~\eqref{eqn:dF_2_BM} for~$\nabla_{(h,h)}^2 \overline{F}[g](t)$.  

Using~\eqref{eqn:equiv_psi} from Theorem~\ref{thm:psi_operator_rep} in the case $n=2$ and~\eqref{eqn:dxQphi} from Lemma~\ref{lemma:heat_gradient_identity}, we obtain (recalling the definition of $v_t(s,x) = (Q_t1)(s,x)$)  
\begin{align}\label{eqn:psi_2}
	&\psi_{t}^{(2)}(t_1,t_2) =  \mbox{$\frac{1}{2}$} (Q_{t_1}\Dx Q_{t_2}\Dx Q_{t}1)(0,0) = \mbox{$\frac{1}{2}$} (Q_{t_1}\Dx Q_{t_2}v_t')(0,0) \nonumber \\
	&= \frac{1}{2} \int_{-\infty}^{g(t_1)}q(0,0;t_1,x_1) (\Dx Q_{t_2}v_t')(t_1,x_1)\,\dif x_1 \nonumber\\
	&=  \frac{1}{2}\int_{-\infty}^{g(t_1)}q(0,0;t_1,x_1) \frac{\partial}{\partial x }\int_{t_2}^t v_t'(u,g(u)-) f_{\tau_g^{t_1,x}}(u)\,\dif u\, \bigg|_{x =x_1}\,\dif x_1\nonumber\\
	&=  \frac{1}{2}\int_{t_2}^t v_t'(u,g(u)-) \int_{-\infty}^{g(t_1)} q(0,0; t_1 ,x_1) \mbox{$\frac{\partial }{\partial x_1  }$}f_{\tau_g^{t_1,x_1}}(u) \,\dif x_1 \,\dif u \nonumber\\
	&= -\frac{1}{2}\int_{t_2}^t v_t'(u,g(u)-) \int_{-\infty}^{g(t_1)}\big(\mbox{$\frac{\partial}{\partial x_1}$} q(0,0; t_1 ,x_1)\big) f_{\tau_g^{t_1,x_1}}(u)\,\dif x_1 \,\dif u \nonumber \\
	&=	- \frac{1}{2}\int_{t_2}^t v_t'(u,g(u)-)K(t_1,u)\,\dif u,
\end{align}
where 
\begin{equation}\label{defn:K_ru}
	K(r, s) :=  \int_{-\infty}^{g(r)}\big(\mbox{$\frac{\partial}{\partial x }$}q(0,0; r ,x )\big) f_{\tau_g^{r,x}}(s) \,\dif x, \quad 0 < r <s.
\end{equation}	
The exchange of order of integrals in the fifth equality is justified using the dominated convergence theorem and the Gaussian bounds for the partial derivatives of the taboo transition density~$q$ (see Lemma 5 in~\cite{Liang2025}). The second last equality was obtained using integration by parts and the boundary conditions~\eqref{eqn:bc_taboo} for the taboo transition density. Using the fundamental theorem of calculus, it follows from~\eqref{eqn:psi_2} that, for $0<t_1<t_2<t,$
\begin{equation}\label{eqn:partial_2_psi}
\mbox{$\frac{\partial}{\partial t_2} $}\psi_{t}^{(2)}(t_1,t_2) =  \mbox{$\frac{1}{2}$}v_t'(t_2,g(t_2)-)K(t_1,t_2).
\end{equation}	
Substituting~\eqref{eqn:partial_2_psi} into~\eqref{eqn:dI} and the resulting expression into~\eqref{eqn:preIBP}, we obtain the claimed expression~\eqref{eqn:dF_2_BM} for~$\nabla_{(h,h)}^2 \overline{F}$.

All that remains is to obtain a closed-form expression for $K(r,u).$ From the Bachelier--L\'evy formula (see e.g.\ 2.0.2 in~\cite[p.~204]{Borodin2015}) we obtain, for $0 < r<s ,$
\begin{equation*}
H(z) := f_{\tau_g^{r,g(r) - z} }(s) =   \frac{z}{s - r}\varphi_{s - r}(z + b(s-r)), \quad z > 0,
\end{equation*}
and
\begin{equation*}
	H'(z) = -\frac{1}{s - r}\bigg(\frac{z^2}{s - r} + b z - 1\bigg)\varphi_{s-r}(z + b(s-r)), \quad z > 0.
\end{equation*}
The taboo transition density of the Brownian motion under the linear boundary~$g$ is given by (see e.g. 1.1.8 on p. 257 in~\cite{Borodin2015}):
\begin{equation}\label{eqn:taboo_transition_density}
	q(s,x;t,y) = \bigl(1 - 
	\exp\bigl\{-\mbox{$\frac{2}{t-s}$}(a+bs -x)(a+bt -y) \bigl\} \bigr)\varphi_{t-s}(y-x), 
\end{equation}
where $x \leq a+bs,$ $y \leq a + bt,$ $0\leq s <t.$ Using the change of variables $z = a + br - x$ and integration by parts we get from~\eqref{defn:K_ru}
\begin{align}\label{eqn:K_step1}
	K(r,s) &=  \int_0^{\infty}\mbox{$\frac{\partial}{\partial y}$}q(0,0; r ,y)|_{y=a+br -z} H(z)\,\dif z \nonumber\\
	&=   \int_0^{\infty}q(0,0; r ,a + br -z) H'(z)\,\dif z \nonumber\\
	&= - \frac{ \gamma}{s - r}\int_0^{\infty}(1-e^{-\frac{2}{r}az})\bigg(\frac{z^2}{s-r} + bz - 1\bigg) e^{-\alpha z^2 + \beta z}\,\dif z, 
\end{align}
where
\begin{equation*}
	\alpha := \frac{s}{2r(s-r)}, \quad \beta := \frac{a}{r}, \quad \gamma := \frac{1}{2\pi \sqrt{r(s-r)}}\exp\bigg\{-\frac{a^2}{2r} - ab - \frac{b^2 }{2} s\bigg\}.
\end{equation*}
For $n=0,1,2,$ let
\begin{equation*}
	M_n(\beta) := \int_0^{\infty}z^n e^{-\alpha z^2 + \beta z}\,\dif z, \quad \beta \in \mathbb{R},
\end{equation*}	
so that we have from~\eqref{eqn:K_step1}, for $0<r<u,$
\begin{equation}\label{eqn:K_step2}
\begin{aligned}
	K(r,u)= -\mbox{$\frac{\gamma}{u-r}$}\Big[&\mbox{$\frac{M_2(\beta) - M_2(-\beta)}{u-r}$}  + b(M_1(\beta) - M_1(-\beta))  - (M_0(\beta) - M_0(-\beta)) \Big].
\end{aligned}
\end{equation}
Using equation 3.462.1 on p.~367 from~\cite{Gradshteyn2015}, the fact that the parabolic cylinder function~$D_{\nu}(z),$ $z \in \mathbb{R},$ $\nu \in \mathbb{Z},$  satisfies $D_{-1}(z) =\sqrt{2\pi } e^{z^2/4}\Phi_1(-z)$ (see e.g.\ equation 9.254.1 on p.\ 1040 in~\cite{Gradshteyn2015} but note that their definition of~$\Phi$ is different to ours), and the recurrence relation $D_{\nu-1}(z)=\nu^{-1}(z D_{\nu}(z) - D_{\nu + 1}(z))$ for negative integers~$\nu$ (see e.g.\ equation 9.247.1 on p.~1039 from~\cite{Gradshteyn2015}),  we obtain
\begin{equation*}
	M_n(\beta) =
	\begin{cases}
		\sqrt{\frac{\pi}{\alpha}}e^{\beta^2/(4\alpha)}\Phi_{2\alpha}(\beta), & n = 0,\\
		\frac{\beta}{2\alpha}M_0(\beta) + \frac{1}{2\alpha}, &n =1,\\
		\frac{\beta^2 + 2\alpha}{4\alpha^2}M_0(\beta)+ \frac{\beta}{4\alpha^2}, & n=2.
	\end{cases}
\end{equation*}
Let $c:= (a\sqrt{2\alpha}/\beta )^{2} = r s /(s-r).$ Substituting
\begin{align*}
	M_0(\beta) - M_0(-\beta) &= \sqrt{\frac{2\pi}{c}}r e^{a^2/(2c)}(2 \Phi_c(a) - 1),\\
	M_1(\beta) - M_1(-\beta) &= a \sqrt{\frac{2 \pi}{c^3}}r^2e^{a^2/(2c)},\\
	M_2(\beta) - M_2(-\beta) &= 2a \frac{r^3}{c^2} + \frac{r^2}{c^2}(a^2 + c)(M_0(\beta) - M_0(-\beta)),
\end{align*}
into~\eqref{eqn:K_step2}, we obtain 
\begin{align}\label{eqn:fla_K_ru_explicit}
	K(r, s) &=  -\frac{a}{s}\varphi_s(a+bs)\bigg[ b + \bigg(\frac{a}{s} - \frac{1}{a} \bigg)(2\Phi_c(a) -1)+  \frac{2c}{s}\varphi_c(a) \bigg] \nonumber\\
		&=-2f_{\tau_g}(s)f(r,s), \quad 0 < r <s,
\end{align}
where $f(r,s)$ was defined in~\eqref{eqn:K(ru)}. This also establishes the asymptotics~\eqref{eqn:K_ru_diagonal} since $c \varphi_c(a) = O(c^{1/2})= O((s-r)^{-1/2}),$ $r\uparrow s.$ Substituting~\eqref{eqn:fla_K_ru_explicit} into~\eqref{eqn:psi_2}, and using~\eqref{defn:first_chaos_linear} for the representation for $\psi_t^{(1)}$, we get
\begin{equation*}
	\psi_t^{(2)}(t_1,t_2) = \int_{t_2}^t f(t_1,s) v_t'(s,g(s)-)f_{\tau_g}(s)\,\dif s = -\int_{t_2}^t f(t_1,s)\mbox{$\frac{\partial}{\partial s}$} \psi_t^{(1)}(s)\,\dif s,
\end{equation*}
and thus the claimed~\eqref{defn:second_chaos_linear} holds.
\end{proof}

\begin{proof}[Proof of Proposition~\ref{prop:approx_m}]
Among the class of linear boundaries $\ell_{a,b}(s) := a + bs,$ $s \geq 0,$ $b \in \mathbb{R},$ with $a = g(0),$ one has $\| g - \ell_{a,b} \|_{\cmSpace_t}^2 = \int_0^t ( \dot{g}(s) - b)^2\,\dif s,$ which is minimised when $b = \frac{1}{t}\int_0^t \dot{g}(s)\,\dif s = \frac{1}{t}(g(t) - g(0)) \equiv b_t.$ Hence $\ell_t \equiv \ell_{g(0),b_t}$ uniquely minimises the $\cmSpace_t$-error and
\begin{equation*}
	\| g - \ell_t\|_{\cmSpace_t}^2 = \| h_t \|_{\cmSpace_t}^2, \quad t >0,
\end{equation*}
where $h_t$ was defined in~\eqref{defn:h_t}. Since $\dot{g}\in\contSpace^1$ by assumption, and $\frac{1}{t}\int_0^t\dot{g}(s)\,\dif s = b_t,$ by Wirtinger's inequality for functions with zero mean (see e.g.\ 258 on p.~185 of~\cite{Hardy1934}) it follows that
\begin{equation*}
	\| g - \ell_t \|_{\cmSpace_t}^2 = \int_0^t ( \dot{g}(s) - b_t)^2\,\dif s \leq \frac{t^2}{\pi^2}\| \ddot{g}\|_{\LpSpace^2((0,t))}^2 .
\end{equation*}
After taking square-roots on both sides, one obtains the claimed inequality for~$\| h_t\|_{\cmSpace_t}.$ Thus~\eqref{ineq:bound_Fm} follows from inequality~\eqref{eqn:error_remainder_taylor} in Theorem~\ref{thm:multi_linear_DF}. 

If $g(t) = a + bt +h(t)$ with $h \in \cmSpace_{\infty},$ then
\begin{equation}\label{bound:dg_b}
	\int_0^t (\dot{g}(s) - b_t)^2 \,\dif s = \int_0^t \dot{h}(s)^2\,\dif s - \frac{h(t)^2}{t} = \| h\|_{\cmSpace_t}^2 - \frac{h(t)^2}{t}.
\end{equation}
By the Cauchy--Bunyakovsky--Schwarz inequality, it follows that $\sup_{s\in[0,t]}\lvert h(s)\rvert \leq t^{1/2}\| h\|_{\cmSpace_t}.$ Hence $h(t)^2/t \leq \|h \|_{\cmSpace_t}^2$ and thus one has from~\eqref{bound:dg_b},
\begin{equation*}
	\sup_{t \in (0,T)}\biggl(\int_0^t (\dot{g}(s) - b_t)^2 \,\dif s\biggr)^{1/2} \leq \| h\|_{\cmSpace_T}.
\end{equation*}
Passing to the limit as $T\to\infty,$ we obtain the final claim.
\end{proof}
 

\section{Examples}\label{sec:examples}
 
In this section we provide five examples which support the theoretical results obtained in Section~\ref{sec:main_result}. 

In Example~\ref{ex:brownian_linear} we numerically assess the sharpness of the claimed bound~\eqref{eqn:error_remainder_taylor} for the truncation error of the Taylor expansion~\eqref{eqn:taylor_series} from Theorem~\ref{thm:multi_linear_DF} in the case of a linear boundary~$g$ and a linear boundary perturbation~$h$. In Example~\ref{ex:fpt_to_chaos} we consider the problem when the first Wiener--It\^o chaos kernel itself is the object of interest and demonstrate how to use~\eqref{eqn:psi_n_fla} from Theorem~\ref{thm:chaos} to obtain numerical approximations for it.  In Example~\ref{ex:OU_chaos_kernel} we use formula~\eqref{eqn:taboo_psi} to numerically obtain the first two chaos kernels in the case of a linear boundary~$g$. 

Next, we compare the accuracy of our proposed approximation $\overline{F}_m$ defined in~\eqref{defn:F_m} with the tangent approximation~\cite{Strassen1967}, Durbin's multiple integral series expansion~\cite{Durbin1992}, and the hazard rate tangent approximation~\cite{Roberts1995}. In Example~\ref{ex:brownian_daniels} we consider the case of the well-known \emph{Daniels boundary}, which has a closed-form expression for the first-passage time distribution. In Example~\ref{ex:brownian_sinusoidal} we turn to the sinusoidal boundary considered in~\cite{Roberts1995} and~\cite{Daniels1996} to demonstrate the effectiveness of our approximation~\eqref{defn:F_m} even when the boundary is not wholly concave nor convex, a case when tangent-type approximations tend to perform worse.

\subsection{Truncation error bounds and chaos kernels}\label{sec:truncation_error_chaos}

\begin{example}[Taylor expansion truncation error bounds]\label{ex:brownian_linear}

In this example we numerically assess the sharpness of the bound~\eqref{eqn:error_remainder_taylor} for the truncation error of the Taylor expansion~\eqref{eqn:taylor_series} from Theorem~\ref{thm:multi_linear_DF} in the case when both the boundary~$g$ and the perturbation~$h$ are linear, a case when we have closed-form formulae.

For $a >0,$ $b \in \mathbb{R},$ let
\begin{equation*}
	g_{a,b}(t) := a +bt, \quad t \in [0,T].
\end{equation*}
It is well-known that the first-passage time distribution tail for the Brownian motion starting at $(s,x),$ $s\in[0,T],$ $x < g_{a,b}(s),$ with the linear boundary~$g_{a,b}$ is given by~\eqref{eqn:v_absx} (see e.g. 1.1.4 on p. 256 in~\cite{Borodin2015}). 

Consider a boundary perturbation $h_c \in \cmSpace$ of the following linear form\textnormal{:}
\begin{equation*}
	h_c(t) = c t, \quad t \in [0,T],\quad c \in \mathbb{R}.
\end{equation*}
Since~$h_c$ is linear, the perturbed boundary $g_{a,b}(t) + h_c(t) = a +(b+c)t$ is also linear, which allows us to compute $\overline{F}[g_{a,b}+h_c](t) = \Px(\tau_{g_{a,b+c}} > t)$ explicitly using~\eqref{eqn:v_absx} by replacing~$b$ with $b+c$. 

Using Theorem~\ref{thm:multi_linear_DF}, it follows that the error arising from approximating $\overline{F}[g_{a,b} + h_c]=\overline{F}[g_{a,b+c}] $ by $\overline{F}[g_{a,b}]$ is given by
\begin{equation}\label{eqn:univariate_remainder}
	 \lvert R_m[g_{a,b},h_c] \rvert= \biggl |\overline{F}[g_{a,b} + h_c] - \sum_{n=0}^{m}\frac{1}{n!}\nabla_{(h_c, \ldots,h_c)}^{n}\overline{F}[g_{a,b}] \biggr|,
\end{equation}
which admits the following upper bound:
\begin{equation}\label{ineq:ex_remainder_bound}
	\lvert R_m[g_{a,b},h_c](t)\rvert \leq \overline{F}[g_{a,b}](t)^{1/2}\sum_{n= m+1}^{\infty}\frac{\lvert c\rvert^nt^{n/2}}{\sqrt{n!}},
\end{equation}
since $\| h_c\|_{\cmSpace_t} = \lvert c\rvert t^{1/2}.$ We can explicitly compute the derivatives $\nabla_{(h_c,\ldots,h_c)}^nF[g_{a,b}]$ in~\eqref{eqn:univariate_remainder} from its definition~\eqref{defn:nabla_hF} using the closed-form expression~\eqref{eqn:v_absx},
\begin{equation*}
	\nabla_{(h_c,\ldots,h_c) }^{n} \overline{F}[g_{a,b}](t) = \frac{\partial^n}{\partial \varepsilon^n} \overline{F}[g_{a,b} + \varepsilon h_c](t)\bigr|_{\varepsilon =0} =
 \frac{\partial^n}{\partial \varepsilon^n}\Px(\tau_{g_{a,b+\varepsilon c}} > t ) \bigr|_{\varepsilon =0},
\end{equation*}
and thus we can compute the remainder~\eqref{eqn:univariate_remainder} explicitly. 

In the case $g_{a,b}(t) = \frac{1}{4} + \frac{3}{10}t$ we plot in  Fig.~\ref{fig:taylor_series} the absolute value of the remainder $R_m[g_{a,b},h_c](1)$ as a function of the number of correction terms~$m$ in the Taylor expansion, together with its upper bound~\eqref{ineq:ex_remainder_bound}. The figure shows that, although the constant factor is somewhat crude, the bound captures the convergence rate of the remainder to zero as $m\to \infty.$ 

\begin{figure}[tbp]
    \centering      
    \includegraphics[scale=0.9]{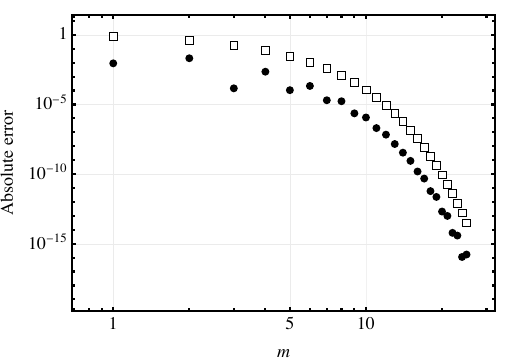}
        \caption{Log-log plot of the absolute value of the Taylor expansion remainder~\eqref{eqn:univariate_remainder} (shown by~$\bullet$) at time $t=1$ and its upper bound~\eqref{ineq:ex_remainder_bound} (shown by~$\square$), where~$m$ is the number of correction terms used in the series expansion~\eqref{eqn:univariate_remainder}. This was done in the case when $a=\frac{1}{4},$ $b=\frac{3}{10}$ and $c = 1.$}
        \label{fig:taylor_series}
\end{figure}

\end{example} 

\begin{example}[Obtaining Wiener--It\^o chaos kernels using first-passage time distributions]\label{ex:fpt_to_chaos}

In this example we will show how to obtain an approximation for the first chaos kernel~$\psi_{t}^{(1)}$ using representation~\eqref{eqn:psi_n_fla} from Theorem~\ref{thm:chaos}. This representation has the advantage that it does not require knowledge of the taboo transition density.

For a small $\varepsilon>0,$ we denote the centered finite difference approximation of the $\cmSpace$-derivative $\nabla_h \overline{F}$ by
\begin{equation}\label{eqn:finite_chaos_1}
	D_h^{\varepsilon}\overline{F}[g]:= \mbox{$\frac{1}{2\varepsilon}$} \big( \overline{F}[g + \varepsilon h] - \overline{F}[g-\varepsilon h]\big) = \nabla_h \overline{F}[g] + O(\varepsilon^2).
\end{equation}
Note that the approximation error in~\eqref{eqn:finite_chaos_1} can be estimated using Theorem~\ref{thm:multi_linear_DF}. 

To numerically evaluate the finite difference~\eqref{eqn:finite_chaos_1} we will further approximate $\overline{F}[g]$ by using the method of integral equations. It is well-known that the first-passage time distribution $F[g]$ for a general $g \in \contSpaceZero$ satisfies the following Volterra integral equation \textnormal{(}see e.g.~\cite{Peskir2002}\textnormal{)}\textnormal{:}
\begin{equation*}\label{eqn:fortet_integral_eqn}
	\overline{\Phi}_t(g(t)) = \int_0^t \overline{\Phi}_{t-s}(g(t) - g(s)) \,\dif F[g](s), \quad t \geq 0,
\end{equation*}
where $\overline{\Phi}_t(x) := 1- \Phi_t(x),$ $ (t,x) \in (0,\infty)\times \mathbb{R},$ and~$\Phi_t$ was defined in~\eqref{defn:varphi}. To obtain an approximation for $F[g],$ one can apply numerical quadrature to the above integral equation and solve the resulting linear system, see e.g.~\cite[p.\ 409]{Peskir2002}. 

The first step to approximately computing the first chaos kernel using representation~\eqref{eqn:psi_n_fla} from Theorem~\ref{thm:chaos} is to compute a vector of finite difference approximations~\eqref{eqn:finite_chaos_1} for boundary perturbation directions~$h(s)$ of the form $\kappa_t(s) = s \wedge t$ for~$t$ taking values on a fixed time grid $\{ t_i  :  0 < t_1 < t_2 < \cdots <t_n < t\}.$ Next, one interpolates this approximation with a continuously differentiable spline, and differentiates the resulting interpolated function to obtain
\begin{equation*}\label{approx:finite_difference}
	\psi_{t}^{(1)}(s) \approx - \frac{\partial}{\partial s} D_{\kappa_s}^{\varepsilon} \overline{F}[g](t),\quad s\in (0,t).
\end{equation*}
Unlike~\eqref{eqn:taboo_psi} which requires knowledge of the taboo transition density and is therefore effectively limited to cases in which that density is available, the finite difference based approximation of~\eqref{eqn:psi_n_fla} described above allows one to obtain approximations for the chaos kernels in the case of general time-dependent boundaries as long as one has a method for accurately numerically computing $\overline{F}[g]$ for arbitrary~$g.$ The tradeoff is that the finite difference based approximation is more computationally expensive. 

In Fig.~\ref{fig:chaos1_compare} we plot the first chaos kernel~$\psi_1^{(1)}$ in the case of the linear boundary $g(t) = \frac{1}{4} + \frac{3}{10} t$ and nonlinear boundary $\tilde{g}(t) := g(t) + \frac{3}{20}\sin(6\pi t).$ This was calculated using the finite difference approximation scheme described above for~\eqref{eqn:psi_n_fla}, where we used $\varepsilon = 10^{-6}$ and a uniform time grid spacing of $t_2 - t_1= 10^{-2}$. 
 
\begin{figure}[!tbp] 
\centering       \includegraphics[scale=0.9]{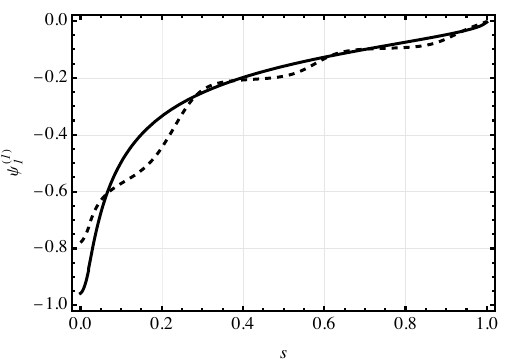}
   \caption{Plot of the first chaos kernel~$\psi_1^{(1)}(s)$ for the cases $g(t) = \frac{1}{4} + \frac{3}{10}t$ (solid black line) and $\tilde{g}(t) = g(t) + 0.15\sin(6\pi t)$ (dashed black line).}
\label{fig:chaos1_compare}
\end{figure}
 \end{example}

\begin{example}[Obtaining Wiener--It\^o chaos kernels using taboo transition densities]\label{ex:OU_chaos_kernel}
In this example we will use representation~\eqref{eqn:taboo_psi} to obtain a numerical approximation for the first two chaos kernels. To do this, one requires a method for computing the taboo transition density~$q$ and the non-crossing probability $(Q_t1)(s,x).$ In the case of a linear boundary~$g(t) = a+bt,$ the taboo transition density and non-crossing probability are given by the closed-form expressions~\eqref{eqn:taboo_transition_density} and~\eqref{eqn:v_absx}, respectively. 

In Fig.~\ref{fig:chaos1} we plot a quadrature approximation of the first chaos kernel 
\begin{equation}\label{eqn:psi_1}
	\psi_t^{(1)}(s) = -\int_{-\infty}^{g(s)}\mbox{$\frac{\partial}{\partial y}$}q(0,0;s,y)(Q_t1)(s,y)\,\dif y, \quad s \in [0,t],
\end{equation}
in the case of a linear boundary $g(t) = a + \frac{3}{10}t$ for vertical intercept values $a = \frac{1}{4},\frac{1}{2},1$. Recalling that the first chaos kernel can be viewed as the `linear sensitivity' of $\mathbf{1}_{\tau_g > t}$ to the increments of the underlying Brownian motion~$W,$ it is probabilistically intuitive that the first chaos kernel is negative since positive increments of~$W$ decrease the likelihood of the path staying entirely below the boundary in the time interval $[0,t]$. We also see from Fig.~\ref{fig:chaos1} that as the boundary recedes away from zero (i.e.\ $a$ becomes larger) the sensitivity `flattens out' near the start of the time interval, which is also intuitive since the path is unlikely to touch the boundary at those earlier times when the boundary is remote. 

\begin{figure}[tbp] 
\centering       
\begin{subfigure}{0.47\linewidth}
\centering       
\includegraphics[scale=0.9]{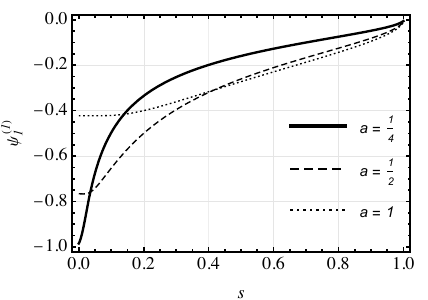}
\caption{First chaos kernel $\psi_t^{(1)}$}
\label{fig:chaos1}
\end{subfigure}\hfill
\begin{subfigure}{0.47\linewidth}
\centering       
\includegraphics[scale=0.9]{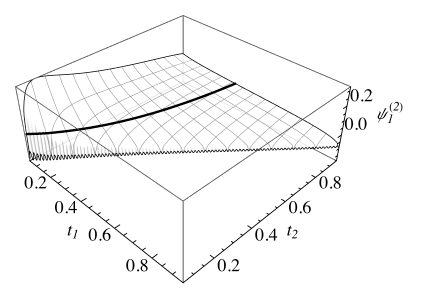}
\caption{Second chaos kernel $\psi_t^{(2)}$}
\label{fig:chaos2}
\end{subfigure}
   \caption{On the left panel we plot the first chaos kernel $\psi_1^{(1)}$ for $g(s) = a + \frac{3}{10}s$ in the cases $a=\frac{1}{4},\frac{1}{2},1$. On the right panel we plot the second chaos kernel~$\psi_1^{(2)}$ for $g(s) = \frac{1}{4} +\frac{3}{10}s.$ For visual reference, the black contour in the plot on the right panel is the level set $\{(t_1,t_2):\psi_1^{(2)}(t_1, t_2) = 0\}.$}
\label{fig:chaos}
\end{figure}

In the case of the second chaos kernel
\begin{equation}\label{eqn:psi_double}
\begin{aligned}
	\psi_t^{(2)}(t_1,t_2)= \frac{1}{2}\int_{-\infty}^{g(t_1)}\int_{-\infty}^{g(t_2)} &( \mbox{$\frac{\partial}{\partial x_1}$}q(0,0;t_1,x_1)) (\mbox{$\frac{\partial}{\partial x_2}$}q(t_1,x_1;t_2,x_2)) \\
	& \times (Q_t1)(t_2,x_2)\,\dif x_2\,\dif x_1, \quad 0 \leq t_1 < t_2<t,
\end{aligned} 
\end{equation}
 we used both~\eqref{eqn:psi_double} and~\eqref{defn:second_chaos_linear} to obtain approximations to good numerical agreement, and plotted the resulting approximation of the latter in Fig.~\ref{fig:chaos2}.  
\end{example}

\subsection{First-passage time distribution approximations}\label{sec:fpt_distr_approx}

In this subsection, we compare the accuracy of our proposed approximation~$\bar{F}_m$ from~\eqref{defn:F_m} with the tangent approximation~\cite{Strassen1967}, Durbin's series expansion~\cite{Durbin1992}, and the hazard rate tangent approximation~\cite{Roberts1995}.

\subsubsection{Tangent approximations}\label{tangent_approx}

For~$g \in \contSpaceZero,$ the classical tangent approximation for the first-passage time density is defined as 
\begin{equation}\label{defn:TA_approx_fpt}
	f^{\text{TA}}(t) := \frac{\lambda(t)}{t}\varphi_t(g(t)), \quad \lambda(t) := g(t) - t \dot{g}(t),\quad t \in (0,T],
\end{equation}
where $\varphi_t$ was defined in~\eqref{defn:varphi}. Justification for the validity of the tangent approximation as~$g$ recedes to infinity for small~$t$ was established in~\cite{Strassen1967} and for general~$t$ in~\cite{Cuzick1981,Jennen1981}. Refinements of this approximation were obtained in~\cite{Jennen1981,Jennen1985}. Error bounds for the approximation were obtained in~\cite{Ferebee1982} and an asymptotic expansion for the error was obtained in~\cite{Ferebee1983}, though its utility remains theoretical rather than numerically practical (see Section~6 in~\cite{Durbin1992}). 

We denote the corresponding approximate distribution function by
\begin{equation}\label{den:F_TA}
	F^{\text{TA}}(t) := \int_0^t f^{\text{TA}}(s) \,\dif s, \quad 	t \in (0,T].
\end{equation}
One of the main drawbacks of the tangent approximation~\eqref{den:F_TA} is that it fails to produce a reasonable approximation when the boundary~$g$ is too steep since the vertical intercept of the tangent approximation can go below the starting value of the Brownian motion. This case is not a problem for our approximation~\eqref{defn:F_m} since our approximating linear boundary~\eqref{defn:g_t} has matching end-point values on $[0,t]$. 

 A series expansion with the tangent approximation being the leading-order term was later obtained and proved to converge pointwise by Durbin~\cite{Durbin1992}. Its uniform convergence was proved in~\cite{Zipkin2011}. Durbin's series~\cite{Durbin1992} is defined as follows: for $m \in \mathbb{N}_{>0},$ let 
\begin{equation}\label{def:durbin_series}
	f_{m}^{\text{Durbin}}(t) := f^{\text{TA}}(t) +\sum_{j=1}^m (-1)^j q_{j+1}(t), \quad  t \in (0,T],
\end{equation}
where, setting $t_0 = t,$ 
\begin{equation*}
\begin{aligned}
	q_j(t) &:= \int_0^t\int_0^{t_1}\cdots \int_0^{t_{j-2}} \frac{\lambda(t_{j-1})}{t_{j-1}} \varphi_{t_{j-1}}(g(t_{j-1}))\\
	&\quad \times 
	\prod_{i = 1}^{j-1}\left[\left(\frac{g(t_{i-1}) - g(t_{i})}{t_{i-1}-t_{i}} - \dot{g}(t_{i-1}) \right)\varphi_{t_{i-1} -t_{i}}(g(t_{i-1}) - g(t_{i}))\right]\dif t_{j-1}\cdots \,\dif t_{2}\,\dif t_1,
	\end{aligned}
\end{equation*}
We denote the corresponding approximate distribution function by
\begin{equation}\label{defn:F_Durbin}
F_{m}^{\text{Durbin}}(t) := \int_0^t f_{m}^{\text{Durbin}}(s)\,\dif s,\quad 	t \in (0,T].
\end{equation}

One of the advantages of Durbin's series is its truncation error $O(c^m/m!),$ which decays to zero faster than our bound $O(c_1^m/\sqrt{m!})$ from~\eqref{eqn:error_remainder_taylor} as the number of terms~$m$ in the series expansion increases. A drawback of Durbin's series is that it is for the first-passage time density and not the distribution function, i.e.\ one needs to integrate term-by-term to obtain the distribution function. This can cause numerical challenges since this increases the number of dimensions of the integrals by one. We note that our approximation~\eqref{defn:F_m} approximates the distribution function directly rather than through the density, reducing the need to integrate one more time.

\subsubsection{Hazard rate tangent approximations}\label{sec:hazard_rate}

The \emph{hazard rate tangent approximation}~\cite{Roberts1995} was introduced as an improvement to the tangent approximation, and is based on approximating the hazard rate of the first-passage time distribution, rather than the first-passage time density. It is defined as
\begin{equation}\label{defn:HRT}
	F^{\text{HRT}}(t) := 1-\exp\left\{ -\int_0^t r(s)\,\dif s \right\},\quad r(t) := \frac{f^{\text{TA}}(t) }{1-F^{\text{TA}}(t)}, \quad 	t \in (0,T].
\end{equation}

One of the major advantages of the hazard rate tangent approximation is its simplicity of implementation. Computationally, the hazard rate tangent approximation compares favourably to both our approximation~\eqref{defn:F_m} and Durbin's series~\eqref{def:durbin_series} since one can obtain an approximation for the distribution function for the entire time interval $t \in [0,T]$ all at once using one-dimensional numerical integration. However, just like the usual tangent approximation, we will see in Example~\ref{ex:brownian_sinusoidal} that the hazard rate tangent approximation fails to produce accurate approximations for the first-passage time density when the boundary is oscillatory.

Other improvements on the tangent approximation include non-linear boundary approximations based on the method of inverse images~\cite{Daniels1996,Lo2002}, which will not be discussed in this paper.
 
\subsubsection{Comparisons}

\begin{example}[Daniels boundary]\label{ex:brownian_daniels}

\begin{figure}[tb] 
\centering       
\begin{subfigure}{0.48\linewidth}
\centering       
\includegraphics[scale=0.9]{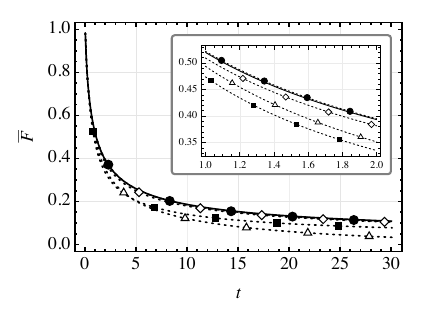}
\caption{Distribution tail}
\label{fig:daniels_approximation}
\end{subfigure}\hfill
\begin{subfigure}{0.48\linewidth}
\centering       
\includegraphics[scale=0.9]{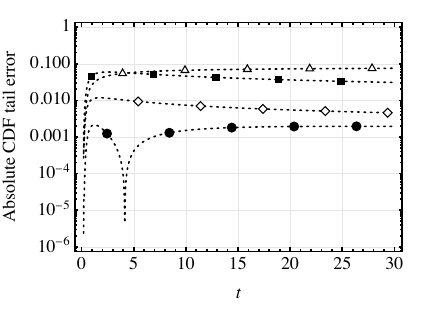}
\caption{Distribution tail error}
\label{fig:daniels_approximation_error}
\end{subfigure}
   \caption{In the left panel we plot the distribution tail~\eqref{eqn:exact_daniels_CDF} of the first-passage time of the Daniels boundary~\eqref{defn:daniels_boundary} (shown by \linestyle[0.8cm]{style5} ), our zeroth-order approximation~$\overline{F}_0$ (shown by~\scalebox{0.65}{$\blacksquare$}), the tangent approximation~\eqref{den:F_TA} (shown by~\scalebox{0.65}{$\triangle$}),  the hazard rate tangent approximation~\eqref{defn:HRT} (shown by~$\diamond$), and our first-order approximation~$\overline{F}_1$ (shown by~$\bullet$) for $t\in [0,30].$ In the left panel we include an inset plot that zooms in on the approximations for $t\in[1,2].$ On the right panel, we have a log-linear plot of the absolute error of the corresponding approximations.}
\label{}
\end{figure}

The so-called \emph{Daniels boundary} was first introduced in~\cite{Daniels1969}, and is defined as
\begin{equation}\label{defn:daniels_boundary}
	g_{\text{D}}(t) := \mbox{$\frac{1}{2}$} - t \ln\Bigl(\mbox{$\frac{1}{4}$}\Bigl(1 + \sqrt{1 + 8e^{-1/t}}\Bigr) \Bigr), \quad t >0,
\end{equation}
where $g_{\text{D}}(0+) =1/2.$ The first-passage time distribution tail for the Daniels boundary~$g_{\text{D}}$ is given by
\begin{equation}\label{eqn:exact_daniels_CDF}
	\overline{F}[g_{\text{D}}](t) = \Phi_t(g_{\text{D}}(t)) - \mbox{$\frac{1}{2}$}\Phi_t(g_{\text{D}}(t) - 1) - \mbox{$\frac{1}{2}$}\Phi_t(g_{\text{D}}(t) -2), \quad t>0,
\end{equation}
see e.g.\ Example 2 on pp.~27--28 in~\cite{Lerche1986}.
We refer the reader to pp.~27--31 in~\cite{Lerche1986} for more examples of time-dependent boundaries which admit closed-form expressions for $\overline{F}[g]$.

In Fig.~\ref{fig:daniels_approximation} we plot the exact distribution tail~\eqref{eqn:exact_daniels_CDF} and the following approximate distribution tails discussed in this paper\textnormal{:} the zeroth and first-order approximations~$\overline{F}_0$ and~$\overline{F}_1$ from~\eqref{defn:F_m}, the tangent approximation~\eqref{den:F_TA}, and the hazard rate tangent approximation~\eqref{defn:HRT}. We did not include the second-order approximation~$\overline{F}_2$ nor Durbin's approximation~\eqref{defn:F_Durbin} in the comparison in Fig.~\ref{fig:daniels_approximation} as it was too computationally expensive to compute for all $t\in[0,30]$ in a reasonable amount of time. The absolute approximation errors of the different methods are shown in Fig.~\ref{fig:daniels_approximation_error}, from which we can see that our approximation~$\overline{F}_1$ is the most accurate in $t\in [0,30].$

In Table~\ref{table:daniels} we report the size of the correction terms in our approximation~\eqref{defn:F_m}, the approximation error, and compare it with the accuracy of Durbin's approximation~\eqref{defn:F_Durbin}. As we mentioned previously, to obtain an approximation for the distribution function using Durbin's approximation one needs to integrate the first-passage time density approximation. Thus a comparison based on the number of correction terms does not fully reflect the computational cost. Comparing the accuracy of our approximation for~$m$ terms with Durbin's approximation for $m-1$ terms, we see that our approximation performs relatively well for $m=1,2.$

\begin{table}[htbp]
\centering
\caption{Results for our approximation~\eqref{defn:F_m} of the boundary crossing probability  $F[g_{\text{D}}](1)= 0.47974935$ up to the displayed number of digits.  The approximation errors for Durbin's series expansion~\eqref{defn:F_Durbin} were obtained from Table~1 in~\cite{Durbin1992} and independently verified by the author.}
\begin{tabular}{crrrrr}\toprule
$m$ & $\frac{1}{m!}\nabla_{(h_1,\ldots,h_1)}^m F[\ell_1](1)$ & $F_m(1)$ & $F_m(1) - F[g_{\text{D}}](1)$ & $F_{m}^{\text{Durbin}}(1) - F[g_{\text{D}}](1)$ \\ \midrule
0 & 
$0.52589901$ & $0.52589901$ & $0.04614966$ & $0.02630264$\\
$1$ & 
$-0.04818438$ & $0.47771463$ & $-0.00203472$ & $-0.00046397$\\
$2$ & 
$0.00183498$ & 
$0.47954962$ & 
$-0.00019974$ &
$-0.00000004$ \\ \bottomrule
\end{tabular}\label{table:daniels}
\end{table}

\end{example}

\begin{example}[Sinusoidal perturbation of a linear boundary]\label{ex:brownian_sinusoidal}

Let
\begin{equation}\label{defn:sin_boundary}
	g(t) := a + bt + \varepsilon \sin(t), \quad t \geq 0, 
\end{equation}
where $a = 2, $ $b = \frac{1}{10}  ,$ $\varepsilon = \frac{1}{4}.$ This oscillating boundary was used in~\cite[(6.2)]{Roberts1995} to test the hazard rate tangent approximation~\eqref{defn:HRT} on a case when the boundary is neither concave nor convex. As noted in~\cite{Roberts1995}, the improvement of the hazard rate approximation~\eqref{defn:HRT} over the tangent approximation~\eqref{den:F_TA} was not theoretically justified in this example. However since $g(t) = a+bt + h(t)$ and $h(t)=\varepsilon\sin(t) \in \cmSpace,$ this example satisfies the conditions in our Propositions~\ref{prop:deriv_bm}--\ref{prop:approx_m}.

\begin{figure}[h] 
\centering       
\begin{subfigure}{0.48\linewidth}
\centering       
\includegraphics[scale=0.9]{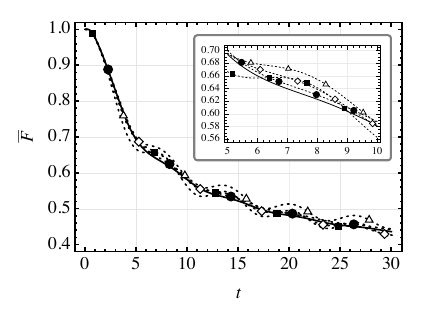}
\caption{Distribution tail}
\label{fig:sin_approximation}
\end{subfigure}\hfill
\begin{subfigure}{0.48\linewidth}
\centering       
\includegraphics[scale=0.9]{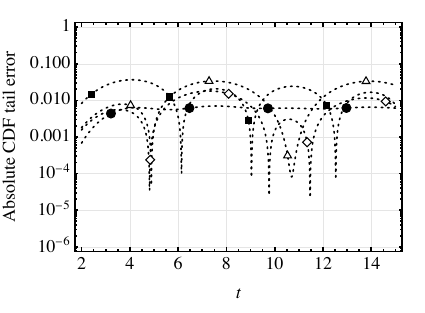}
\caption{Distribution tail error}
\label{fig:sin_approximation_error}
\end{subfigure}
   \caption{In the left panel we have a plot of the distribution tail of the first-passage time of the sinusoidal boundary~\eqref{defn:sin_boundary} (shown by \linestyle[0.8cm]{style5} ), our zeroth-order approximation~$\overline{F}_0$ (shown by~\scalebox{0.65}{$\blacksquare$}), the tangent approximation~\eqref{den:F_TA} (shown by~\scalebox{0.65}{$\triangle$}),  the hazard rate tangent approximation~\eqref{defn:HRT} (shown by~$\diamond$), and our first-order approximation $\overline{F}_1$ (shown by~$\bullet$) for $t\in [0,30].$ In the left panel we include an inset plot of the approximations for $t\in[5,10].$ On the right panel, we have a log-linear plot of the absolute error of the corresponding approximations.}
\label{}
\end{figure} 

Unsurprisingly,  due to the oscillations in the boundary~\eqref{defn:sin_boundary}, we can see from Fig.~\ref{fig:sin_approximation} that the tangent approximation~\eqref{den:F_TA} fares poorly in this example. The hazard rate tangent approximation~\eqref{defn:HRT} works relatively better, however we see from Fig.~\ref{fig:sin_approximation_error} that the approximation error oscillates. As observed in~\cite{Roberts1995}, these successive periods of under- and overestimation cancel out. This, however, suggests that the corresponding first-passage time density approximation may be less accurate.

We approximate the first-passage-time density by taking numerical derivatives with respect to~$t$ for our approximation~\eqref{defn:F_m}, and denote it by~$f_{1}(t) := -\frac{\dif}{\dif t}\overline{F}_{1}(t),$ $t >0.$ The hazard rate tangent approximation for the first-passage time density is given by
\begin{equation}\label{defn:fpt_HRT}
	f^{\text{HRT}}(t) := r(t)(1 - F^{\text{HRT}}(t) ), \quad t \in (0,T).
\end{equation}
In Fig.~\ref{fig:sin_fpt_approx} we plot the first-passage time density approximations $f_1,$ $f_1^{\text{Durbin}},$ $f^{\text{HRT}},$ $f^{\text{TA}}.$ From Fig.~\ref{fig:sin_fpt_approximation_error} we see that our approximation~$f_1(t)$ has the smallest absolute error among the four approximations for sufficiently large~$t$ in this particular example. These promising results indicate further research on the first-passage time density approximation~$f_1$ may be fruitful.

\begin{figure}[!tbp] 
\centering       
\begin{subfigure}{0.48\linewidth}
\centering       
\includegraphics[scale=0.9]{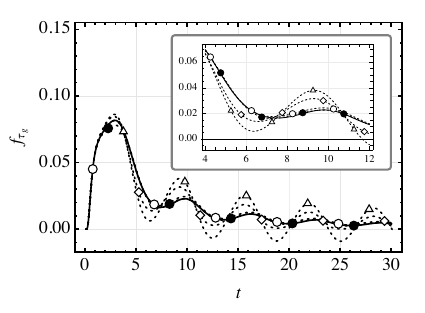}
\caption{First-passage time density}
\label{fig:sin_fpt_approx}
\end{subfigure}\hfill
\begin{subfigure}{0.48\linewidth}
\centering       
\includegraphics[scale=0.9]{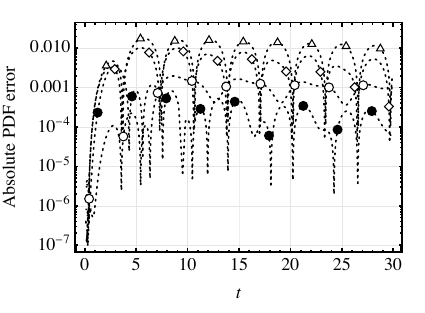}
\caption{First-passage time density error}
\label{fig:sin_fpt_approximation_error}
\end{subfigure}
   \caption{In the left panel we have a plot of the first-passage time density for the sinusoidal boundary~\eqref{defn:sin_boundary} (shown by \linestyle[0.8cm]{style5} ), Durbin's approximation~$f_1^{\text{Durbin}}$ from~\eqref{def:durbin_series} (shown by~$\circ$), the tangent approximation $f^{\text{TA}}$ from~\eqref{defn:TA_approx_fpt} (shown by~\scalebox{0.65}{$\triangle$}),  the hazard rate tangent approximation~$f^{\text{HRT}}$ from~\eqref{defn:fpt_HRT} (shown by~$\diamond$), and our first-order approximation~$f_1$ (shown by~$\bullet$) for $t\in [0,30].$ In the left panel we include an inset plot of the approximations for $t\in[4,12].$ On the right panel, we have a log-linear plot of the absolute error of the corresponding approximations.}
\label{}
\end{figure} 
\end{example}

\textbf{Acknowledgments:} The author gratefully thanks Persi Diaconis for scanning and providing a copy of M.\ Hogan's PhD thesis. The author also thanks Kostya Borovkov for numerous helpful discussions and suggestions related to boundary crossing problems, Roxanne He and Panpan Ren for their editorial advice on the draft, and Neil Chada and Laurent Mertz for supporting the author's postdoctoral appointment.

\end{document}